\documentclass{amsart}

\usepackage{microtype}
\usepackage{pinlabel}  %%% the recommended graphics+labelling package
\usepackage{graphicx}  %%% the recommended graphics package
\usepackage{amscd}
\usepackage{amsmath}
\usepackage{bm}
\usepackage{enumerate}
\usepackage[dvipdfmx]{hyperref}

\title[Simple closed geodesics on hyperelliptic translation surfaces]
{Simple closed geodesics on hyperelliptic translation surfaces 
and classification theorem for translation surfaces in $\calhhyp(4)$}

\author{Yoshihiko Shinomiya}
\address{Mathematics Education, 
Faculty of Education College of Education, 
Academic Institute, Shizuoka University 
836 Ohya, Suruga-ku, Shizuoka 422-8529, JAPAN}
\email{shinomiya.yoshihiko@shizuoka.ac.jp}
\urladdr{}


\newtheorem{theorem}{Theorem}[section]    % Standard theorem environment

\newtheorem{proposition}[theorem]{Proposition}   

\newtheorem{lemma}[theorem]{Lemma}          % Lemma environment with numbering 
\theoremstyle{definition}
\newtheorem{definition}[theorem]{Definition}    % Definition environment with 
\newtheorem{remark}{Remark}             % Unnumbered environment for remarks.
\newtheorem{example}[theorem]{Example}   

\theoremstyle{plain}
 \newtheorem{corollary}[theorem]{Corollary}

\newtheorem*{mthm1}{Theorem \ref{main}}

\newtheorem*{mthm3}{Theorem \ref{classification}}
\newtheorem*{mcor4}{Corollary \ref{classification_cor}}

\newtheorem*{ackn}{Acknowledgements}
\newcommand{\thm}[1]{\begin{theorem} #1\end{theorem}}
\newcommand{\prop}[1]{\begin{proposition} #1\end{proposition}}
\newcommand{\lem}[1]{\begin{lemma} #1\end{lemma}}
\newcommand{\rem}[1]{\begin{remark} #1\end{remark}}
\newcommand{\pf}[1]{\begin{proof} #1\end{proof}}
\newcommand{\defi}[1]{\begin{definition} #1\end{definition}}
\newcommand{\exam}[1]{\begin{example} #1\end{example}}
\newcommand{\cor}[1]{\begin{corollary} #1\end{corollary}}

\newcommand{\ord}{\mathrm{ord}}
\newcommand{\calh}{\mathcal{H} }
\newcommand{\calhhyp}{\mathcal{H}^{\mathrm{hyp}}}
\newcommand{\calheven}{\mathcal{H}^{\mathrm{even}}}
\newcommand{\calhodd}{\mathcal{H}^{\mathrm{odd}}}
\newcommand{\calhnonhyp}{\mathcal{H}^{\mathrm{nonhyp}}}

\newcommand{\vv}{\mathbf{v}}
\newcommand{\interior}[1]{{\mathrm{Int}}(#1)}

\newcommand{\sltz}{\mathrm{SL}(2, \mathbb{Z})}
\newcommand{\sltr}{\mathrm{SL}(2, \mathbb{R})}
\newcommand{\gltr}{\mathrm{GL}(2, \mathbb{R})}
\newcommand{\str}{{\mathrm{St}}_5}
\newcommand{\veech}{\Gamma(X,\omega)}
\newcommand{\ccyl}{C_{\mbox cyl}}
\newcommand{\teich}{Teichm\"{u}ller }
\begin{document}

\begin{abstract}    % type your abstract below
In this paper, we give the maximum of 
the numbers $n$ 
such that we can take $n$  
simple closed geodesics without singularities 
that are  disjoint to each other 
for translation surfaces in the hyperelliptic components 
$\calhhyp(2g-2)$ and $\calhhyp(g-1, g-1)$.
The maximum is different from that of the case of hyperbolic surfaces.
%of disjoint and non-homotopic 
%closed geodesics which do not contain singularities 
%for translation surfaces in the hyperelliptic components 
%$\calhhyp(2g-2)$ and $\calhhyp(g-1, g-1)$. 
%%%%%%%%%%%%%%%%%%%%%%%%%%%%%%%%%%%%%%%%%%%%%%%%%%%%%%%%%
%We  show that translation surfaces in the  hyperelliptic components 
%are always described from certain polygons constructed from parallelograms. 
We also give a classification theorem for translation surfaces in $\calhhyp(4)$ 
with respect to their Euclidean structures.
\end{abstract}

\maketitle

%%%%%%%%%%%%%%%%%%%%   Start of main body of article

\section{Introduction}
A translation surface $(X, \omega)$ is a Riemann surface $X$ 
together with a holomorphic $1$-form $\omega$ on $X$.
The $1$-form $\omega$ induces a singular Euclidean structure on $X$. 
We can also construct a Euclidean structure on a compact Riemann surface $X$ 
from a holomorphic quadratic differential $q$ on $X$. 
The pair $(X, q)$ is called a flat surface.  
Translation surfaces and flat surfaces are important  for \teich theory. 
The \teich space ${\mathcal T}_g$ is the space of equivalence classes 
of marked Riemann surfaces of genus $g$.  
The \teich space ${\mathcal T}_g$  has a complete metric called \teich metric. 
By the uniformization theorem, \teich space ${\mathcal T}_g$ is also regarded 
as the space of equivalence classes of marked hyperbolic structures on a surface of genus $g$. 
From this, \teich spaces are studied through the hyperbolic geometry on surfaces. 
%For example, we can induce the Fenchel–Nielsen coordinates on \teich spaces.
%On a hyperbolic surface of genus $g\geq 2$, every homotopy class 
%of simple closed curves contains
%a unique geodesic. 
%There always exists a $3g-3$ 
%which is a universal covering of the moduli space.
However,  deformations of Riemann surfaces along geodesics on the \teich space ${\mathcal T}_g$ 
are difficult to understand from hyperbolic structures. 
The \teich theorem claims that \teich geodesics on ${\mathcal T}_g$ 
are represented by the deformation of translation surfaces and flat surfaces.

The moduli space $\calh_g$ of translation surfaces of genus $g$ 
is stratified by the orders of 
zeros of $1$-forms. 
Fix $g \geq 2$. 
For positive integers $k_1, \dots, k_m$ satisfying 
$k_1 \leq \cdots \leq k_m$ and 
$k_1+ \cdots + k_m=2g-2$, 
the stratum $\calh(k_1, \dots, k_m)$ 
consists of all equivalence classes of pairs $(X, \omega)$ 
such that $\omega$ has exactly $m$ zeros of orders $k_1, \dots, k_m$.
Kontsevich-Zorich \cite{KonZor03} proved that  $\calh (k_1, \dots, k_m) $ 
has at most $3$ connected components. 
The stratum $\calh(2g-2)$ and $\calh(g-1, g-1)$ have connected components
called hyperelliptic components. 
We denote them by $\calhhyp(2g-2)$ and $\calhhyp(g-1,g-1)$, respectively. 
A translation surface $(X, \omega) \in \calh (2g-2)$ (resp. $\calh (g-1, g-1)$)  
is an element of $\calhhyp(2g-2)$ (resp. $\calhhyp (g-1, g-1)$)
if and only if $X$ has the hyperelliptic involution $\tau$ and it satisfies $\tau^\ast \omega=-\omega$.  
We call translation surfaces in  $\calhhyp (2g-2) \cup\calhhyp (g-1, g-1)$ \emph{hyperelliptic translation surfaces}.
In the case of genus $2$, we have $\calh(2)=\calhhyp(2)$ and $\calh(1, 1)=\calhhyp(1,1)$. 
%AHyperelliptic translation surfaces of genera $2$ are investigated by 

In this paper, our main interests are hyperelliptic translation surfaces and their simple closed geodesics.
On a translation surface, a closed geodesic may pass through  some singularities. 
We call simple closed geodesics without singularities \emph{regular}.  
One of the main problems in this paper is 
to determine the maximum of the numbers $n\in \mathbb{N}$
such that there exist 
regular closed geodesics $\gamma_1, \dots, \gamma_n$ 
which are disjoint and not homotopic  
to each other on every hyperelliptic translation surface $(X, \omega)$. 
We denote the maximum number by $N(X, \omega)$.
In the case of hyperbolic surfaces of genus $g\geq 2$, 
the maximum number of pairwise disjoint simple closed geodesics is always $3g-3$. 
It is given by pants curves 
and the maximum number depends only on the genus $g$. 
However, the maximal numbers depend on 
their strata for translation surfaces.
Nguyen \cite[Theorem 2.8]{Nguyen17} investigated the case of genus $2$. 
If $(X, \omega)\in \calh(2)$, then we have $N(X, \omega)=2$.
If $(X, \omega)\in \calh(1, 1)$, then we have $N(X, \omega)=3$.
One of our main theorems is a generalization of this result to general genus.
\begin{mthm1}
Let $g\geq 2$.
If $(X, \omega) \in \calhhyp(2g-2)$, then we have $N(X, \omega)=g$.
If $(X, \omega) \in \calhhyp(g-1, g-1)$, then we have $N(X, \omega)=g+1$. 
\end{mthm1}
We can easily show by a topological argument that 
$N(X, \omega) \leq g$ for $(X, \omega) \in \calhhyp(2g-2)$ 
and 
$N(X, \omega) \leq g+1$ for $(X, \omega) \in \calhhyp(g-1, g-1)$.
To prove the converse inequality, 
we give $g$ or $g+1$ regular closed geodesics 
which are disjoint to each other.
This is done by constructing a good triangulation of 
the quotient flat surface $(Y, q)$ of $(X, \omega)$ 
by the hyperelliptic involution.
The flat surface $(Y, q)$ is of genus $0$ 
and $q$ is a holomorphic quadratic differential 
on $Y$ 
that has only one zero and some simple poles.
The triangles are  ``restricted triangles'' defined as in Definition \ref{restricted_triangle}.
All vertices of this triangulation are the unique zero of $q$ and 
all poles of $q$ are midpoints of some edges of triangles.
This triangulation gives us 
saddle connections 
such that they connect poles and they are disjoint to each other. 
%The ``good'' triangulation 
%has vertices only on the zero of $q$ and 
%gives us $g$ or $g+1$ 
%segments on $(Y, q)$ which connect poles and are disjoint to each other.
The preimages of the saddle connections are regular closed geodesics which we want. 
By investigating 
triangulations of $(X, \omega)\in  \calhhyp(4)$
by  ``restricted triangles'',  
we obtain the following theorem.
% for $\calhhyp(4)$.
Here, $\str$ is a translation surface in $\calhhyp(4)$ constructed from five squares as in Figure \ref{str} 
and $\gltr \cdot \str$ is its $\gltr$-orbit. 
A simple cylinder is a Euclidean cylinder in a translation surface that is foliated by parallel regular closed geodesics 
and each of whose boundary components passes the singular point only once.
\begin{mthm3}
Let $(X, \omega) \in \calhhyp(4)$ 
and $\tau$ the hyperelliptic  involution of $(X, \omega)$.
If $(X, \omega) \not \in \gltr \cdot \str$, 
then there exist disjoint simple cylinders $C_1, C_2, C_3$ each of which is 
invariant under $\tau $. 
If $(X, \omega) \in \gltr \cdot \str$, 
then $(X, \omega)$ has at most two  disjoint simple cylinders. 
\end{mthm3}
%This theorem gives us other constructions of  translation surfaces in  $\calhhyp(4)$. 
The following corollary is a classification theorem for translation surfaces in  $\calhhyp(4)$
with respect to the Euclidean structure on translation surfaces.
%As the classification theorem for topological surfaces, 
Every translation surface  $(X, \omega) \in \calhhyp(4) \setminus \gltr \cdot \str$ 
and its  Euclidean structure
are realized by gluing three Euclidean cylinders with a flat surface of genus $0$.

\begin{mcor4}[Classification Theorem for translation surfaces in $\calhhyp(4)$]
 Let $(X, \omega) \in \calhhyp(4) \setminus \gltr \cdot \str$. 
Then, there exists a flat surface $(X_0, q)$ of genus $0$ 
such that  
 $q$ 
has a zero $p_0$ of order $2$ and six simple poles $p_1, p_2, \dots, p_6$, 
an involution $\tau_0 : X_0 \to X_0$, 
saddle connections $s_1, s_2, s_3$, and 
Euclidean cylinders $C_1, C_2, C_3$
satisfying the following (see Figure \ref{classification3}):
\begin{enumerate}
\item $\tau_0^\ast q=q$, 
\item $s_i$ is a saddle connection of $(X_0, q)$ connecting $p_0$ and $p_i$ for $i=1, 2, 3$, 
\item $\tau_0(s_i)$ connects $p_0$ and $p_{i+3}$ for $i=1, 2, 3$, 
\item two of each $s_1, s_2, s_3, \tau_0(s_1), \tau_0(s_2)$, and $\tau_0(s_3)$ 
         intersect  only at $p_0$,
\item the circumference of $C_i$ equals $2|s_i|$ for $i=1,2,3$, and
\item  $(X, \omega)$ is obtained by cutting $(X_0, q)$ along the saddle connections 
$s_1$, $s_2$, $s_3$, $\tau_0(s_1)$, $\tau_0(s_2)$, $\tau_0(s_3)$ 
and gluing $C_i$ to the slits $s_i$ and $\tau_0(s_i)$ for all $i=1, 2, 3$.
\end{enumerate}
\end{mcor4}

\begin{ackn}
This work was supported by JSPS KAKENHI Grant Number 17K14184.
The author would like to thank Professor Yoshihide Okumura 
for leading and fostering the author to a mathematician.
The author also would like to thank Hidetoshi Masai 
for his valuable comments.
\end{ackn}

\section{Preliminaries}
Through this paper, 
we assume that Riemann surfaces are always 
analytically finite.
For an analytically finite Riemann surface $X$, we denote by $\overline{X}$
the compact Riemann surface $X$ with the punctures filled.

\subsection{Holomorphic $1$-forms and their moduli $\calh_g$}
\defi{[Holomorphic $1$-forms] 
Let $X$ be a Riemann surface.
A holomorphic $1$-form $\omega$ on $X$ is 
a tensor whose restriction  to every local coordinate neighborhood $(U, z)$
has the form $f dz$, 
where $f$ is a holomorphic function on $U$. 
}

Let $\omega$ be a holomorphic $1$-form 
on a Riemann surface $X$.
Choose local coordinate neighborhoods $(U, z), (V, w)$ 
with $U \cap V \not =\emptyset$.
If $\omega=f dz$ on $(U, z)$ and 
$\omega=g dw$ on $(V, w)$, we have
$$f (z) \left(\frac{dz}{dw}\right)=g(w).$$
Thus, if $p_0\in U\cap V$ is a zero of $f$ of order $n$ 
then it is also a zero of $g$ of order $n$.
Now, zeros of $\omega$ and their orders are well-defined.
We denote by ${\rm Sing}(\omega)$ the set of all zeros of $\omega$. 
Let  $\ord_\omega (p)$ denote the order of a zero $p$ of $\omega$. 
A holomorphic $1$-form on $X$ whose restriction to 
every  local coordinate neighborhood  
is  constant function $0$ is denoted by $0$.
By the Riemann-Roch formula, we have the following.
\prop{\label{RHf}
Let $X$ be a compact  Riemann surface of genus $g \geq 2$ and
$\omega \not =0$ a holomorphic $1$-form on $X$. 
Then we have 
$$\sum_{p\in {\rm Sing}(\omega)} \ord_\omega  (p) =2g-2.$$
}

Fix $g\geq 2$.
Let ${\mathcal M}_g$ be the moduli space of
compact Riemann surfaces of genus $g$. 
For a compact Riemann surface $X$,  
we denote by $\Omega(X)$
the vector space of all holomorphic $1$-forms on $X$.  
We set $\Omega^\ast (X)=\Omega(X)-\left\{0\right\}$.
Then, the moduli space of pairs $(X, \omega)$ is defined by
$$\calh_g=\left\{ (X, \omega) : X\in {\mathcal M}_g, \omega \in \Omega^\ast (X) \right\}/\sim.$$
Where, for pairs $(X, \omega)$ and $(Y, \omega^\prime)$, 
the relation $(X, \omega) \sim (Y, \omega^\prime)$ holds
if there exists a conformal map $f: X \rightarrow Y$ such that 
$f^\ast \omega^\prime =\omega$.
The moduli space $\calh_g$ is a complex algebraic orbifold 
of dimension $4g-3$. 
This is fibered over the moduli space ${\mathcal M}_g$ with
the fiber over each $[X] \in {\mathcal M}_g$ equals $ \Omega^\ast (X) / {\rm Aut} (X)$. 

The orbifold $\calh_g$ is stratified by the orders of 
zeros of $1$-forms. 
Let $k_1, \dots, k_m$ be positive integers whose sum is $2g-2$ 
and which satisfy $k_1 \leq k_2 \leq \dots \leq k_m$.
The subspace $\calh(k_1, \dots, k_m)$ of $\calh_g$ 
consists of all equivalence classes of pairs $(X, \omega)$
such that $\omega$ has exactly $m$ zeros of orders $k_1, \dots, k_m$. 
Then we have 
$$\calh_g =\bigsqcup_
{\substack{ 0<k_1 \leq k_2 \leq \dots \leq k_m \\
k_1+\cdots + k_m=2g-2}  }
\calh(k_1, \dots, k_m).$$
Each  $\calh(k_1, \dots, k_m)$ is a stratum of $\calh_g$.
It is known that the dimension of  
 $\calh(k_1, \dots, k_m)$ 
is $2g+m-1$ 
(see  \cite{Masur82}, \cite{Veech82} and \cite{Veech90}).

\exam{
The orbifold $\calh_2$ has two strata $\calh(2)$, and $\calh (1,1)$. 
The orbifold $\calh_3$ has five strata 
$\calh(4), \calh(2, 2), \calh(1, 3), \calh(1, 1, 2)$, and $\calh(1, 1, 1, 1) $.
}

In the rest of this section, we study connected components of strata of $\calh_g$.
In this paper, our interests are hyperelliptic components.

\defi{
Let $g \geq 2$.
Let $\calhhyp (2g-2)$ be the subset of $\calh(2g-2)$ consisting of all  $(X, \omega)$ 
such that $X$ has the hyperelliptic involution $\tau$ 
and $\omega$ satisfies $\tau^\ast \omega=-\omega$. 
Let $\calhhyp (g-1, g-1)$ be the subset of $\calh(g-1, g-1)$ consisting of all $(X, \omega)$
such that $X$ has the hyperelliptic involution $\tau$ 
and $\omega$ satisfies $\tau^\ast \omega=-\omega$. 
}

\rem{
For the case of genus $2$, 
we have 
$\calh(2)=\calhhyp(2)$ 
and 
$\calh(1, 1)=\calhhyp(1, 1)$.
}

The connected components  of the strata are completely classified 
by Kontsevich and Zorich. 
See \cite{KonZor03} for the definition of  
$\calheven, \calhodd$, and 
$\calhnonhyp$.

\thm{
[\cite{KonZor03}, Theorem 2]
In the case of genus two, two strata $\calh(2)$ and $\calh(1, 1)$ are connected. 
In the case of genus three,  
$\calh (4)$ has two connected components
 $\calhhyp (4)$ and  $\calhodd(4)$. 
 The stratum $\calh(2,2)$ has two  connected components
$\calhhyp (2, 2)$ and $\calhodd(2, 2)$.
The other strata are connected.
}

 \thm{
 [\cite{KonZor03}, Theorem 1] 
Let $g \geq 4$. 
We have the following. 
\begin{itemize}
\item The stratum $\calh (2g-2)$ has three  connected components 
 $\calhhyp (2g-2), 
 \calheven(2g-2)$, and $\calhodd(2g-2)$.
 
\item If $g$ is odd,  the stratum $\calh(g-1,g-1)$ has three  connected components 
$\calhhyp (g-1, g-1 ),  \calheven(g-1, g-1)$, and $\calhodd(g-1, g-1)$.

\item If $g$ is even,  the stratum $\calh(g-1, g-1)$ has two connected components
$\calhhyp (g-1, g-1)$ and $\calhnonhyp(g-1, g-1)$. 

\item All other strata which is of the form $\calh(2k_1, \dots, 2k_n)$ 
$(k_1, \dots, k_n \in \mathbb{N})$ 
has two  connected components 
 $\calheven(2k_1, \dots, 2k_n)$ and $\calhodd(2k_1, \dots, 2k_n)$.

\item The others are connected.
\end{itemize}
 }

 For $g \geq 2$, we call $\calhhyp (2g-2)$ and $\calhhyp(g-1, g-1)$ \emph{hyperelliptic components}.
 If $(X, \omega) \in \calhhyp (2g-2)$,  the unique zero of $\omega$ 
 is fixed by the hyperelliptic involution of $X$. 
%Thus, it is a Weierstrass point of $X$.
If $(X, \omega) \in \calhhyp (g-1, g-1)$, two zeros of $\omega$ are 
 permuted by the hyperelliptic involution of $X$.
%Thus, they are not Weierstrass points of $X$.

\subsection{Meromorphic quadratic differentials}
\defi{[Meromorphic quadratic differential] 
Let $X$ be a Riemann surface.
A meromorphic quadratic differential $q$ on $X$ is 
a tensor whose restriction to every local coordinate neighborhood $(U, z)$
has the form $f dz^2$, where $f$ is a meromorphic function on $U$. 
For a meromorphic quadratic differential $q$ on $X$,  
we denote by $|q|$ 
the differential $2$-form on $X$ whose restriction to 
 every local coordinate neighborhood $(U, z)$ 
 is $|f|dx dy$ if  the restriction of $q$  to $(U, z)$ is $f dz^2$. 
Here, $x$ is the real part of $z$ and $y$ is the imaginary part of $z$.
The $L^1$-norm $||q||$ of a meromorphic quadratic differential $q$ on $X$ is 
defined by $||q||=\iint_X |q|$.
}

We assume that meromorphic quadratic differentials  
always have finite $L^1$-norms.
By this assumption, meromorphic quadratic differential $q$ 
may have poles of order at most $1$. 
% on the punctures of $X$. 
Thus, the poles of $q$ are simple poles if they exist.
%For convenience, we sometimes call  poles of $q$ of orders $1$  zeros of orders $-1$ if they exist. 
The set of all zeros and poles of $q$ is denoted by ${\rm Sing}(q)$. 
%Note that if $q$ has a pole, ${\rm Sing}(q)$ is not a subset of $X$ but of $\overline{X}$.

Let $X$ be a compact Riemann surface of genus $g$.
Let $q$ be a meromorphic quadratic differential on a Riemann surface $X$.
The order of a zero $p$ of $q$ is denoted by $\ord_q (p)$.
If $p$ is a simple pole of $q$, we set $\ord_q(p)=-1$. 
For other point $p^\prime$ of $X$, we set $\ord_q(p^\prime)=0$. 
By the Riemann-Roch formula,  we have
$$\sum_{p\in X} \ord_q (p) =4g-4.$$
%Here, $n$ is the number of poles of $q$.

\subsection{Translation surfaces and Flat surfaces}

%Here, we define flat surfaces and translation surfaces. 
%We give two equivalent ways to define them. 
%One is to use a holomorphic quadratic differentials or a holomorphic $1$-forms. 
%The other is to use an atlas.

\defi{[Translation surface and flat surface]
A translation surface $(X, \omega)$ is a pair of a Riemann surface $X$ 
and a non-zero holomorphic $1$-form $\omega$ on $X$. 
The points of ${\rm Sing}(\omega)$ are called singular points of $(X, \omega)$.
A flat surface $(X, q)$ is a pair of a Riemann surface $X$ 
and a non-zero meromorphic quadratic differential $q$ on $X$.
The points of ${\rm Sing}(q)$ are called singular points of $(X, q)$.
}

Let $(X, \omega)$ be a flat surface.
The holomorphic $1$-form $\omega$ gives a Euclidean structure on $X$.
If $p_0\in X$ is not a zero of $\omega$, we may choose a neighborhood $U$ of $p_0$
such that the map
$$\displaystyle z(p)=\int_{p_0}^p \omega :  U\to \mathbb{C}$$
is a chart of $X$. 
The collection of such charts $u=\left\{ (U, z)\right\}$
is an atlas on $X-{\rm Sing}(\omega)$.
The transition functions of this atlas are always of the forms $w= z+({\rm const.})$.
Thus, $u$ is a singular Euclidean structure on $(X, \omega)$.
If $p_0$ is a zero of $\omega$ of order $n$ ($n=1, 2, \dots$), 
% a zero $p_0 \in {\rm Sing}(q)$, 
there exists a chart $(U, \zeta)$ of $X$ around $p_0$ 
such that $\omega$ is represented as $\omega=\zeta^n d\zeta$.
Then, we have a chart 
$$\displaystyle 
z(p)=\int_{p_0}^p \omega
=\int_{p_0}^p \zeta^n  d\zeta
=\frac{1}{n+1} \zeta (p) ^{n+1}
$$
around $p_0$.
With respect to this chart, the angle around $p_0$ is $2(n+1)\pi$. 
%The points in ${\rm Sing}(q)$ are called singular points of $(X, q)$.

Let $(X, q)$ be a translation surface. 
An atlas $u$ on $X-{\rm Sing}(q)$ is also given by $q$. 
For each $p_0 \in X-{\rm Sing}(q)$, 
we can choose a neighborhood $U$ such that the integration 
$$\displaystyle z(p)=\int_{p_0}^p \sqrt{q} :  U\to \mathbb{C}$$
is well-defined and gives  a chart.
The collection of such charts $u=\left\{ (U, z)\right\}$ 
is an atlas on  $X-{\rm Sing}(q)$  whose 
transition functions  are of the forms 
$w=\pm z+({\rm const.})$. 
If $p_0$ is a zero of $q$ or a pole of $q$ with $\ord_q(p)=n$ ($n=-1, 1, 2, \dots$), 
% a zero $p_0 \in {\rm Sing}(q)$, 
there exists a chart $(U, \zeta)$ of $X$ around $p_0$ 
such that $q$ is represented as $q=\zeta^n d\zeta^2$.  
Then, we have a chart 
$$\displaystyle 
z(p)=\int_{p_0}^p \sqrt{q}
=\int_{p_0}^p \zeta^\frac{n}{2}   d\zeta
=\frac{2}{n+2} \zeta (p) ^\frac{n+2}{2}
$$
around $p_0$.
This implies that the angle around $p_0$ with respect to this chart is $(n+2)\pi$.

%Replacing $\omega$ of the above integrations to $\sqrt{q}$, 
%the $1$-form $\omega$ gives an atlas on $X$ whose transition functions
%Hence, $\omega$ also gives a singular Euclidean structure on $X$. 
%The points in ${\rm Sing}(\omega)$ are also called singular points of $(X, \omega)$.

\rem{ 
By the definition of $\calh_g$, 
the space $\calh_g$ is also regarded  as the moduli space 
of translation surfaces of genus $g$.
}

%Next, we give another definition of translation surfaces and flat surfaces. 
%
%\defi{[Translation surface and flat surface]
%A translation surface (resp. flat surface) $(X, u)$ is a pair of a Riemann surface $X$ 
%and an atlas on $X-Z$ for some finite subset of $X$ satisfying the following: 
%\begin{enumerate}
%\item $u$ is compatible to the conformal structure of $X$,  
%\item the transition functions are always of the form 
%\begin{align*}
% w=z+({\rm const}.) &
%\end{align*}
% (resp. $w=\pm z+ ({\rm const}.$) )  
% for 
%$(U, z), (V, w)\in u$ with $U \cap V \not = \emptyset$,
%\item $u$ is the maximal with respect to the above two conditions.
%\end{enumerate}
%The atlas $u$ is called a translation (resp. flat) structure on $X$.
%}
%Let $(X, u)$ be a translation surface such that 
%$u$ is a chart on $X-Z$. 
%By definition, we have a  holomorphic $1$-form 
%$\omega$ whose restriction to a chart $(U, z)\in u$ 
%is $\omega=dz$. 
%The maximality  of $u$ implies that  $Z={\rm Sing}(\omega)$.   
%Thus above two definitions are equivalent.
%

\subsection{Regular closed geodesics and curve complexes}
We define some terminologies for translation surfaces 
which are also defined for flat surfaces as well. 

\defi{[Saddle connection]
A saddle connection  
is a geodesic segment on  a translation surface $(X, \omega)$ 
joining two singular points or a singular point to itself 
and containing no singular points in its interior.
}

For any hyperbolic surface, 
every homotopy class of simple closed curves 
contains a unique geodesic. 
This is not true for translation surfaces and flat surfaces.

\prop{[\cite{Strebel84}, Theorem 14.3]\label{geod_rep}
Let $(X, \omega)$ be a translation surface. 
Let $\gamma$ be a simple closed curve on $(X, \omega)$.
Then, one of the following holds. 
\begin{itemize}
\item The homotopy class of $\gamma$ 
contains a unique closed geodesic of $(X, \omega)$.  
Moreover, the geodesic is a concatenation of saddle connections. 
The angles between consecutive
saddle connections are at least $\pi$ 
on both sides.

\item The homotopy class of $\gamma$ 
contains infinitely many simple closed geodesics containing no singular points.
They are parallel to each other and
the union of all such geodesics forms an open Euclidean cylinder $C_\gamma$.
Each of the boundary components  of  $C_\gamma$ is  
a concatenation of saddle connections 
that are parallel to $\gamma$.
\end{itemize}  
}

%Let $\gamma$ be a simple closed geodesic  on $(X, q)$
%which contains no singular points. 
%By 
%there always exists the Euclidean cylinder $C_\gamma$.

\defi{[Regular closed geodesic]
A simple closed geodesic $\gamma$  on a translation surface 
% on $(X, q)$
that contains no singular points
is called regular. 
%For a regular closed geodesic $\gamma$, 
%the Euclidean cylinder $C_\gamma$ as in Proposition \ref{geod_rep} 
%is called the geometric cylinder for $\gamma$. 
}

%\begin{prob}
% Let $(X, \omega)$ be a translation surface.  
%We set 
%\begin{align*}
%N(X, \omega)
%=\max \left\{ n \in \mathbb{N} : 
%\exists \mbox{regular closed geodesics }  \gamma_1, \dots, \gamma_n 
%\mbox{ s.t. disjoint and not homotopic  
%to each other}
%\right\} 
%\end{align*} 
%We denote by $N(X, \omega)$ 
%the maximum of numbers $n\in \mathbb{N}$
%such that there exist 
%regular closed geodesics $\gamma_1, \dots, \gamma_n$  
%on $(X, \omega)$
%which are disjoint and not homotopic  
%to each other. 
%\end{prob}

Let $(X, \omega)$ be a translation surface. 
We denote by $N(X, \omega)$ 
the maximum of the numbers $n\in \mathbb{N}$
such that there exist 
regular closed geodesics $\gamma_1, \dots, \gamma_n$  
on $(X, \omega)$
which are disjoint and not homotopic  
to each other. 
%the maximum number of 
%disjoint regular closed geodesics 
%which are not homotopic to each other. 
For a hyperbolic surface of genus  $g$ with $n$ punctures, 
the maximum number of disjoint closed geodesics is always 
$3g-3+n$. 
Hence, the maximum numbers depend only on the topological types of surfaces. 
However, for translation surfaces, it is not true.

\thm{[\cite{Nguyen17}, Theorem 2.8]\label{Nguyen_thm}
Let $(X, \omega)$ be a translation surface of genus $2$.
If $(X, \omega) \in \calh(2)$, then $N(X, \omega)=2$. 
If $(X, \omega) \in \calh(1, 1)$, then $N(X, \omega)=3$. 
}

Theorem \ref{Nguyen_thm} is also  
described  in terms of curve complex in \cite{Nguyen17}.
The curve complex $C(S)$ for  a surface $S$ of finite type is a simplicial complex such that 
the vertices are the free homotopy classes of essential simple closed curves on $S$ 
and $k+1$ vertices span a simplex if and only if  they have pairwise disjoint representatives 
for each  $k \geq 0$. 
The curve complexes for translation surfaces are analogies of this.

\defi{[Curve complex for translation surface $(X, \omega)$]
The curve complex $\ccyl=\ccyl (X, \omega)$ for a translation surface $(X, \omega)$ 
is a simplicial complex such that  
the vertices are regular closed geodesics on $(X, \omega)$ 
and $k+1$ vertices span a simplex if and only if  they are  pairwise disjoint for each  $k \geq 0$. 
}

\thm{[\cite{Nguyen17}, Theorem A]\label{Nguyen_thm2}
Let $(X, \omega)$ be a translation surface of genus $2$.
If $(X, \omega) \in \calh(2)$, 
then $\ccyl(X, \omega)$ contains no $2$-simplexes.
If $(X, \omega) \in \calh(1, 1)$, 
then $\ccyl(X, \omega)$ contains some $2$-simplexes.
}

One of our main results  is a generalization of these theorems 
to the case of higher genus case.

\thm{\label{main}
Let $g\geq 2$.
If $(X, \omega) \in \calhhyp(2g-2)$, then $N(X, \omega)=g$. 
If $(X, \omega) \in \calhhyp(g-1, g-1)$, then $N(X, \omega)=g+1$. 
} 

\cor{
Let $g\geq 2$.
If $(X, \omega) \in \calhhyp(2g-2)$, 
then $\ccyl(X, \omega)$ contains some $(g-1)$-simplexes
but contains no $g$-simplexes.
If $(X, \omega) \in \calhhyp(g-1, g-1)$, 
then $\ccyl(X, \omega)$ contains some $g$-simplexes 
but contains no $(g+1)$-simplexes.
}

\subsection{$\sltr$-orbits of translation surfaces and Veech groups}
Let $g \geq 2$.
There exists an action of $\sltr$ on $\calh_g$ which leaves each stratum invariant. 
Let $A\in \sltr$ and $(X, \omega) \in \calh_g$. 
Then, $\omega$ induces an Euclidean structure $u=\left\{ (U, z)\right\}$ 
on $X$. 
We regard $A$ as a linear map and compose it to each chart of $u$.
Then, the atlas $A\circ u =\left\{ (U, A\circ z)\right\}$ 
is a Euclidean structure on  the surface $X$.
We consider it as a new conformal structure on $X$.
Then, we have a holomorphic $1$-form $A\circ \omega$ on $X$ 
whose restriction to each chart of  $(U, w) \in A\circ u$ is $dw$. 
The $1$-forms $\omega$ and $A\circ \omega$ are in the same stratum.
The $\sltr$-orbit of  a translation surface $(X, \omega)$ is defined by
\begin{align*}
\sltr  \cdot (X, \omega)=\left\{(X, A\circ \omega) : A\in \gltr\right\}. 
\end{align*}

\defi{[Veech group]
The Veech group $\veech$ of $(X, \omega)$ is the subgroup of all $A \in \sltr$ 
which leaves $(X, \omega)$ invariant.
That is, $\veech$ is written as 
\begin{align*}
\veech =\left\{ A\in \sltr : (X, A\circ \omega) =(X, \omega) \right\}.
\end{align*}
}

%If $(X, \omega)$ has finite Euclidean area

\thm{[\cite{Veech89}, Proposition 2.7
%, see also \cite{HerSch07}
]
The Veech group $\veech$ of a translation surface $(X, \omega)$ 
is a discrete subgroup of $\sltr$. 
}

Veech groups are important to understand translation surfaces.

\defi{[periodic direction]
Let $(X, \omega)$ be a translation surface.
A direction $\theta \in [0, \pi)$ is called periodic 
if  every $z\in X$ is contained 
in a regular closed geodesic of direction $\theta$ or a saddle connection  
of direction $\theta$.
}

\thm{[Veech's dichotomy Theorem, \cite{Veech89}, 
Propositions 2.4, 2.10, and 2.11
]\label{veech_dic}
Let $(X, \omega)$ be a translation surface 
whose Veech group $\veech$ is a lattice in $\sltr$
 (i.e. a subgroup of finite covolume). 
For each direction $\theta \in [0, \pi)$, one of the following holds:
\begin{enumerate}
\item\label{periodic} $\theta$ is periodic, 
\item every geodesic with direction $\theta$ is dense in $X$.
\end{enumerate}
Moreover, in case (\ref{periodic}), there exists a parabolic element $A\in \veech$ 
which has a unique eigenvector of direction $\theta$. 
If $\theta=0$, the fixed point of $A$ as a  M\"{o}bius transformation is $\infty$. 
If $\theta \neq 0$, the fixed point of $A$ as a  M\"{o}bius transformation is $\cot \theta$.
}

We denote by ${\rm Peri}(X, \omega)$ the subset of $[0, \pi)$ consisting of all periodic directions of $(X, \omega)$.
Let us define the action of $\sltr$ on $[0, \pi)$.
For any $A\in \sltr$ and $\theta, \theta^\prime \in [0, \pi)$, 
we assume that 
the equation $A\cdot \theta=\theta^\prime$ holds 
if and  only if 
\begin{align*}
A \left[
\begin{array}{c} 
\cos \theta\\
\sin \theta 
\end{array}
\right]
=t
\left[
\begin{array}{c} 
\cos \theta^\prime\\
\sin \theta^\prime 
\end{array}
\right]
\end{align*}
holds 
for some $t \in \mathbb{R}$.
With respect to this action, 
the Veech group $\veech$
preserves 
 ${\rm Peri}(X, \omega)$.
Moreover, we have the following by Theorem \ref{veech_dic}.

\prop{\label{peripara}
Let $f: [0, \pi)
 \to \hat{\mathbb{R}}$ be a function such that 
\begin{align*}
f(\theta)=
\begin{cases}
\cot \theta & (\theta \in (0, \pi)) \\
\infty &  (\theta=0).
\end{cases}
\end{align*} 
Let $(X, \omega)$ be a translation surface 
whose Veech group $\veech$ is a lattice in $\sltr$.
Then, the set $f\left({\rm Peri}(X, \omega) \right)$ coincides with 
the set of all fixed points of parabolic elements of $\veech$. 
Moreover, there exists a bijection from 
${\rm Peri}(X, \omega)/\Gamma(X, \omega)$ 
to $f({\rm Peri}(X, \omega))/\Gamma(X, \omega)$. 
Here $\Gamma(X, \omega)$ acts on ${\rm Peri}(X, \omega)$ 
as linear maps and 
on $f({\rm Peri}(X, \omega))$ as M\"{o}bius transformations.
} 

\pf{
Let  $\theta, \theta^\prime \in {\rm Peri}(X, \omega)$. 
Assume that there exists 
$A=\left[
\begin{array}{ccc} 
a & b \\ 
c & d 
\end{array}
\right]\in \Gamma(X, \omega)$ 
satisfying $A \cdot \theta =\theta^\prime$.
Then the equation  
\begin{align*}
A \left[
\begin{array}{c} 
\cos \theta\\
\sin \theta 
\end{array}
\right]
=t
\left[
\begin{array}{c} 
\cos \theta^\prime\\
\sin \theta^\prime 
\end{array}
\right]
\end{align*}
holds for some $t \in \mathbb{R}$.
This implies the equation 
\begin{align*}
 \frac{a f(\theta)+b}{c f(\theta)+d}=f(\theta^\prime)=f(A\cdot \theta)
\end{align*}
holds. 
Therefore, the function $f$ induces a surjective map from 
${\rm Peri}(X, \omega)/\Gamma(X, \omega)$
to 
$f({\rm Peri}(X, \omega))/\Gamma(X, \omega)$. 
By following this argument in reverse, we can also conclude that this map is injective.

Next, we show that 
the set $f\left({\rm Peri}(X, \omega) \right)$ coincides with 
the set of all fixed points of parabolic elements of $\veech$. 
By Theorem \ref{veech_dic},  there exists 
a parabolic element 
$
A=\left[
\begin{array}{ccc} 
a & b \\ 
c & d 
\end{array}
\right]
\in \veech$ 
which has a unique eigenvector of direction $\theta$. 
By the same way as above
\begin{align*}
\frac{a f(\theta)+b}{c f(\theta)+d}=f(\theta) 
\end{align*}
holds.
Therefore, $f(\theta)$ is a fixed point of a parabolic element $A \in \veech$. 
The converse is also proved by following this argument in reverse.
}

Let $(X, \omega)$ be a translation surface 
whose Veech group $\veech$ is a lattice in $\sltr$.
Assume that  $\theta \in [0, \pi)$ is periodic.
By Theorem \ref{veech_dic},  
the union of all regular closed geodesics of direction $\theta$ 
consists of some open Euclidean cylinders $C_1, \dots, C_n$. 
Each boundary component of the cylinders consists of saddle connections of direction $\theta$.
The family of the cylinders $\left\{ C_1, \dots, C_n\right\}$ is called 
the \emph{cylinder decomposition} for $\theta$.

%
%Let $(X, \omega)$ be a translation surface 
%whose Veech group $\veech$ is a lattice in $\sltr$.

\prop{\label{affine_homeo}
Let $(X, \omega)$ be a translation surface and 
$u=\left\{(U, z)\right\}$ the atlas of $(X, \omega)$.
Let $\left\{C_1, \dots, C_m \right\}$ and $\left\{ C_1^\prime, \dots, C_n^\prime \right\}$ 
be cylinder decompositions for periodic directions  
$\theta$ and $\theta^\prime$ of $(X, \omega)$, respectively.
Suppose that there exists $A\in \veech$ which maps the direction $\theta$ to $\theta^\prime$. 
Then, we have $m=n$ and 
 there exists a homeomorphism $h: (X, \omega) \to (X, \omega)$ 
%  and 
%a permutation $\sigma \in S_n$ 
satisfying the following:
 \begin{enumerate}
\item $\left\{h(C_1), \dots, h(C_m)\right\} =\left\{ C_1^\prime, \dots, C_n^\prime \right\}$ 
%$h(C_i)=C_{\sigma(i)}^\prime$  for all $i=1, 2, \dots, n$ 
and
\item\label{affine} for any $(U, z)$ and $(V, w) \in u$ with $h(U)\subset V$, 
the composition $w \circ h \circ z^{-1}$ is an affine map 
whose derivative is $A$.
\end{enumerate}
%which is an affine map with respect to the atlas $u$ induced by $\omega$.

%$A\left(\left\{ C_1, \dots, C_m\right\} \right)=\left\{ C_1^\prime, \dots, C_n^\prime \right\}$. 
%Here, $A$ is a map between cylinders whose universal covering is the linear map $A: \mathbb{C} \to \mathbb{C}$.
}

\pf{
Let $H_i$ be the height of the cylinder $C_i$ and 
$\widetilde{C_i}$ the $\left( \frac{1}{2}H_i \right)$-neighborhood  
of the line 
$\left(\sin \theta  \right) x - \left( \cos \theta \right)y=0$ in $\mathbb{C}$ 
for each $i=1, 2, \dots, m$.
Then $\widetilde{C_i}$ is the universal covering of $C_i$.
The linear map $\widetilde{h_i} : \widetilde{C_i}  \to \widetilde{h_i} (\widetilde{C_i}) ; z \mapsto A z$ 
is projected to a homeomorphism $h_i : C_i \to h_i(C_i)$.
Gluing the cylinders $h_1(C_1), \dots, h_m(C_m)$ in the same way 
as the construction of $(X, \omega)$ from $C_1, \dots, C_m$, we obtain a flat surface $(X, A\circ \omega)$. 
Since $A\in \Gamma(X, \omega)$, we have $(X, A\circ \omega)=(X,  \omega)$.
This implies that $\left\{h_1(C_1), h_2(C_2), \dots, h_m(C_m)\right\}$ coincides with the cylinder decomposition 
$\left\{ C_1^\prime, \dots, C_n^\prime \right\}$.
%Since $C_i$ is foliated by regular closed geodesics of direction $\theta$, 
%$h_i(C_i)$ is foliated by regular closed geodesics of direction $\theta^\prime$.
%By assumption that $\theta^\prime \in {\rm Peri}(X, \omega)$, 
%the family $\left\{h_1(C_1), h_2(C_2), \dots, h_m(C_m)\right\}$ coincides with the cylinder decomposition 
%for $\theta^\prime$. 
Thus, we have $m=n$.
Now, we can define a homeomorphism $h: (X, \omega) \to (X, \omega)$ 
such that $h|_{C_i}=h_i$ for all $i=1, 2, \dots, n$ 
and satisfies the condition (\ref{affine}). 
%Hence, there exists  $\sigma \in S_n$ satisfying 
%$h(C_i)=C_{\sigma(i)}^\prime$ for all $i=1, 2, \dots, n$.
}

\begin{figure}[ht!]
\labellist
%\small
\hair 0pt
\pinlabel $1$  at 195 160
\pinlabel $1$  at  50 160
\pinlabel $2$  at -10 95
\pinlabel $2$  at  135 95
\pinlabel $3$  at -10 32
\pinlabel $3$  at  75 32
\pinlabel $a$ at 155 195
\pinlabel $a$ at 155 110
\pinlabel $b$ at 95 50
\pinlabel $b$ at 95 195
\pinlabel $c$ at 35 135
\pinlabel $c$ at 35 -10
\endlabellist
\centering
 \includegraphics[scale=0.4]{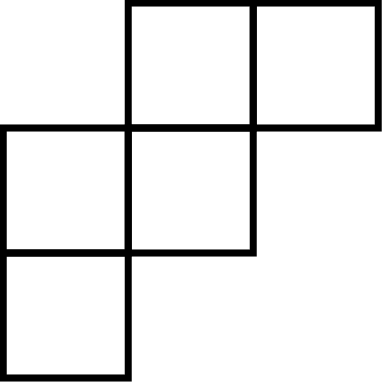} 
\caption{The translation surface $\str$.}
\label{str}
\end{figure}

We finish this section after 
 giving an example of Veech groups.
\defi{
[\cite{Schmithusen06}]
The translation surface $\str \in\calhhyp(4)$ is
constructed from $5$ squares as in Figure \ref{str}  and the $1$-form $dz$ on it 
with the same labels glued.
}

\rem{
The action of the hyperelliptic involution $\tau$ on $\str$  
is seen by cutting $\str$ into five squares 
and rotating them by $\pi$ rotation about the centers. 
We can glue the resulting squares in the same way as $\str$. 
Then we again have $\str$. 
The centers of the squares correspond to the fixed points of $\tau$.
The other fixed points are the unique singular point and 
the points corresponding to the midpoints of sides labeled by $a$ and $3$ 
as in Figure \ref{str}.
}

\exam{\label{exam_st5}
The Veech group $\Gamma(\str)$ is  calculated 
by Schmith\"{u}sen (\cite[Proposition 13]{Schmithusen06}). 
%She use an algorithm for calucutlating Veech gruops  \cite{Schmithusen04}). 
The Veech group $\Gamma(\str)$ is described as 
\begin{align*}
 \Gamma(\str)
 =\left\{ 
\left[
\begin{array}{ccc} 
a & b \\ 
c & d
\end{array}
\right]
\in \sltz : 
a+c \mbox{ and } 
b+d \mbox{ are odd } 
 \right\}.
\end{align*}
Set 
$
T= \left[
\begin{array}{ccc} 
1 & 1\\
0 & 1
\end{array}
\right]
$
and 
$
R=
\left[
\begin{array}{ccc} 
0 & -1\\
1 & 0
\end{array}
\right]
$.
By using the Reidemeister-Schreier method (see \cite{Schmithusen04}), 
we can see that 
\begin{align*}
 \Gamma(\str)
 =\left<
 R, T^2, (TR)T(TR)^{-1}
 \right>
\end{align*}
and the coset representatives of $\Gamma(\str)$ in $\sltr$ 
is $\left\{I, T, T R \right\}$. 
The quotient $\mathbb{H}/\Gamma(\str)$ 
is an orbifold  of genus $0$ and has two cusps and only one cone point of order $2$.
We observe the action of homeomorphisms corresponding to 
generators $R$, $T^2$, and $(TR)T(TR)^{-1}$ 
as in Proposition \ref{affine_homeo}. 
Let $P$ be the union of $5$ squares.
The action of $R$ is seen by rotating $P$  
by $\frac{\pi}{2}$-rotation and cutting it into five squares and pasting them 
so that the form is the same as $P$ (see Figure \ref{St_5R}).
To see the action of $T^2$, we cut $\str$ along the horizontal sides of squares.
Then $\str$ is decomposed into $3$ cylinders. 
On the cylinder of area $1$, the matrix $T^2$ acts 
as the square of the right-hand Dehn twist along the core curve. 
The action of $T^2$ onto the other cylinders 
is the right-hand Dehn twists along the core curves (see Figure \ref{St_5T}). 
The actions induce an action of $T^2$ on $\str$ 
since each boundary component of the cylinders is pointwise fixed by the actions.  
The action of $(TR)T(TR)^{-1}$ is difficult to see. 
Let $S\in \sltr$ be a $\frac{\pi}{4}$-rotation. 
By computation, we have $(TR)T(TR)^{-1}=ST^2S^{-1}$.
This implies that $(TR)T(TR)^{-1}$ preserves the direction $\frac{\pi}{4}$.
Let us cut $P$ along the diagonals of direction  $\frac{\pi}{4}$ 
and glue them as in Figure \ref{St_5another}.
Let $P^\prime$ be the resulting parallelogram.
Regard the bottom vertex of  $P^\prime$ as the origin.
We act the linear map  $(TR)T(TR)^{-1}=ST^2S^{-1}$ to $P^\prime$.
Then, the permutation of the labels of the upper side  given by $(TR)T(TR)^{-1}$ is 
a cyclic permutation 
$\sigma=
 \begin{pmatrix}
 e & d & c & b & a
 \end{pmatrix}
$. 
The labels of the lower side of $P^\prime$ in invariant under $(TR)T(TR)^{-1}$. 
Next, by identifying the sides with label $f$, the translation 
$z \mapsto 
z+2(1+i) 
$
on $P^\prime$ 
corresponds to 
the permutation 
$\tau=
 \begin{pmatrix}
 a & d & b & e & c
 \end{pmatrix}
$
for the labels of the upper side and 
the permutation 
$\mu=
 \begin{pmatrix}
 a & c & e & b & d
 \end{pmatrix}
$
for the labels of the lower side.
Therefore, the composition of the two actions 
corresponds to 
the permutation 
$\tau \sigma =
 \begin{pmatrix}
 a & c & e & b & d
 \end{pmatrix}
$
for the labels of the upper side and 
the permutation 
$\mu=
 \begin{pmatrix}
 a & c & e & b & d
 \end{pmatrix}
$
for the labels of the lower side. 
As $\tau \sigma=\mu$, 
 the composition of the two actions is well-defined on the translation surface $\str$.
\begin{figure}[ht!]
\labellist
%\small
\hair 0pt
\pinlabel $1$ at 30 30
\pinlabel $2$ at 30 90
\pinlabel $3$ at 90 90
\pinlabel $4$ at 90 150
\pinlabel $5$ at 150 150
\pinlabel $R$ at 245 110
\pinlabel $\rotatebox{90}{1}$ at 455 30
\pinlabel $\rotatebox{90}{2}$ at 395 30
\pinlabel $\rotatebox{90}{3}$ at 395 90
\pinlabel $\rotatebox{90}{4}$ at 335 90
\pinlabel $\rotatebox{90}{5}$ at 335 150
\pinlabel $=$ at 520 90
\pinlabel $\rotatebox{90}{5}$ at 595 30
\pinlabel $\rotatebox{90}{4}$ at 595 90
\pinlabel $\rotatebox{90}{3}$ at 655 90
\pinlabel $\rotatebox{90}{2}$ at 655 150 
\pinlabel $\rotatebox{90}{1}$ at 715 150
\endlabellist
\centering
 \includegraphics[scale=0.4]{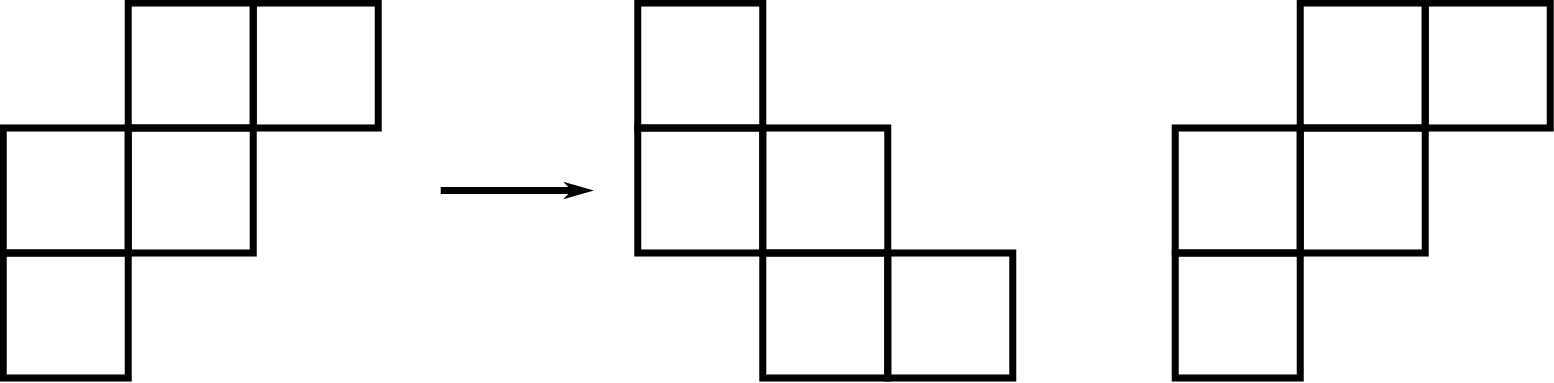} 
\caption{The action of $R$ onto $\str$.}
\label{St_5R}
\end{figure}

\begin{figure}[ht!]
\labellist
%\small
\hair 0pt
\pinlabel $T$ at 250 125
\endlabellist
\centering
 \includegraphics[scale=0.4]{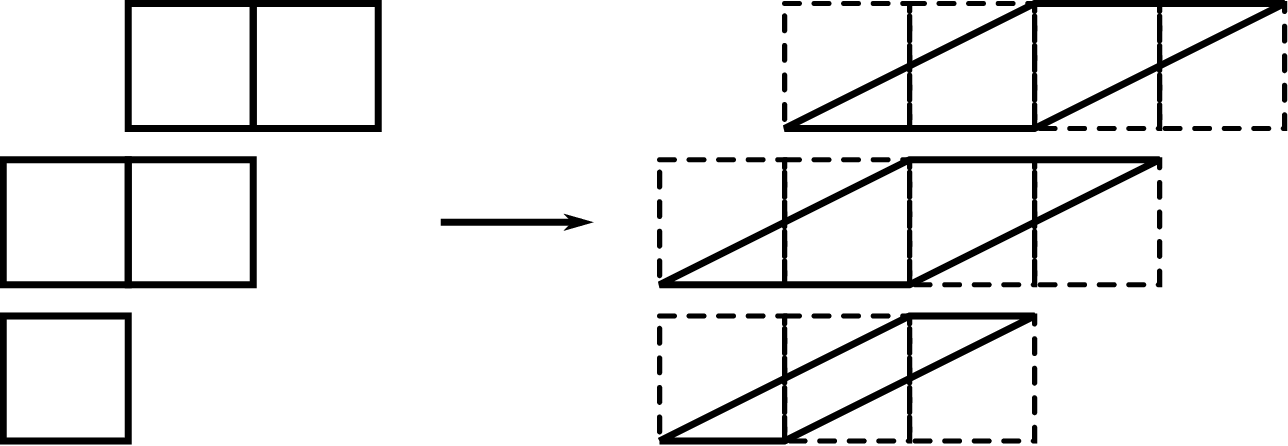} 
\caption{The action of $T$ onto $\str$.}
\label{St_5T}
\end{figure}

\begin{figure}[ht!]
\labellist
%\small
\hair 0pt
\pinlabel $1$ at 15 145
\pinlabel $1^\prime$ at 45 130
\pinlabel $2$ at 15 208
\pinlabel $2^\prime$ at  45 187
\pinlabel $3$ at 77 208
\pinlabel $3^\prime$ at 107 187
\pinlabel $4$ at 77 268
\pinlabel $4^\prime$ at  107 247
\pinlabel $5$ at 139 268
\pinlabel $5^\prime$ at 169 247
\pinlabel $1$ at 320 45
\pinlabel $2^\prime$ at 345 80
\pinlabel $3$ at 380 105
\pinlabel $4^\prime$ at 405 140
\pinlabel $5$ at 440 165
\pinlabel $5^\prime$ at 470 205
\pinlabel $4$ at  500 225
\pinlabel $3^\prime$ at 525 265
\pinlabel $2$ at 560 285 
\pinlabel $1^\prime$ at  585 320  
\pinlabel $\mbox{cut and paste}$ at 285 225
\pinlabel $e$ at 320 102
\pinlabel $d$ at 383 162
\pinlabel $c$ at 446 222
\pinlabel $b$ at 500 282
\pinlabel $a$ at 562 342
\pinlabel $a$ at 345 22
\pinlabel $b$ at 408 82
\pinlabel $c$ at 471 142
\pinlabel $d$ at 525 202
\pinlabel $e$ at 585 262

\pinlabel $f$ at 287 30
\pinlabel $f$ at 625 335
\endlabellist
\centering
 \includegraphics[scale=0.4]{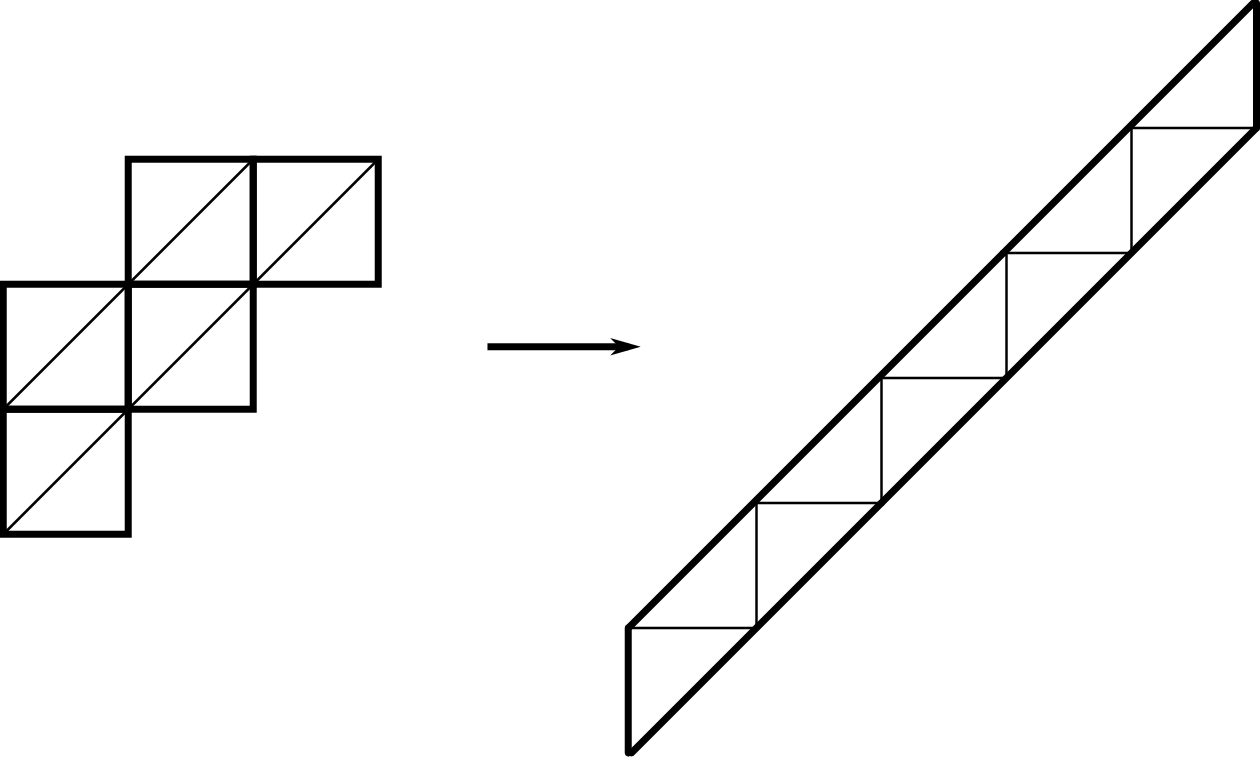} 
\caption{Another description of $\str$ with labels $a, b, c, d, e$, and $f$.} 
\label{St_5another}
\end{figure}
}
Observing the actions of elements of $\Gamma(\str)$, we have the following which is used 
for Lemma \ref{only_two}.

\prop{\label{preserve_fp}
Let $\tau$ be the hyperelliptic involution of $\str$. 
Let $p_1$ and $p_2$ be fixed points of $\tau$ which are 
the images of the midpoints of sides with labels $a$ and $3$, respectively. 
For any $A\in \Gamma(\str)$, the action of $A$ onto $\str$ 
preserves the set $\left\{p_1, p_2\right\}$. 
}

\pf{ We can easily see that 
the generators $R, T^2, (TR)T(TR)^{-1}$ of $\Gamma(\str)$ 
preserve the set $\left\{p_1, p_2\right\}$ 
by observing the above example. 
Moreover, $p_1$ and $p_2$ are fixed points of $\tau$.
Thus, we have the claim.
}
%\prop{ 
%Let $(Y, q)$ be a flat surface. 
%Let $s$ be a segment from a singular point $p_0$
%of $(Y, q)$ to a simple pole $p_1$ of $q$.
%and $\vv$ the vector corresponding to $s$.
%
%
%}

\section{Proof of Theorem \ref{main}}\label{proof_of_main_thm}
In this section, we prove Theorem  \ref{main}.
% and  Theorem \ref{net_thm}. 
Let  $g\geq 2$.
Let  $(X, \omega)$ be a translation surface in $\calhhyp (2g-2)\cup \calhhyp(g-1, g-1)$
with the hyperelliptic involution $\tau$.  

\subsection{Upper estimates of $N(X, \omega)$}
\prop{\label{upper_bound_proof}
If $(X, \omega)\in \calhhyp(2g-2)$, then we have $N(X, \omega) \leq g$. 
If $(X, \omega)\in \calhhyp(g-1,g-1)$, then we have $N(X, \omega) \leq g+1$. 
}

\pf{
Assume that $(X, \omega)\in \calhhyp(2g-2)$. 
If $N(X, \omega) \geq g+1$, then there exist  disjoint regular closed geodesics 
$\gamma_1, \dots, \gamma_{g+1}$ in $(X, \omega)$ that are not homotopic to each other.
Since $(X, \omega)$ is a surface of genus $g$, 
$X-\left(\gamma_1 \cup \cdots \cup \gamma_{g+1}\right)$ 
has at least 
two connected components.
Since $(X, \omega)$ has only one singular point, one of the connected components $X^\prime$ 
is a translation surface 
that has 
%only one geodesic boundary and 
no singular points. 
%Such translation surfaces are only simply connected 
Taking a double of $X^\prime$, we obtain a compact translation surface without singular points.
Then, it must be a torus.
This implies that $X^\prime$ is a Euclidean cylinder.
Now, the boundary components of $X^\prime$ are regular closed geodesics 
that are homotopic to each other.
%However, $X^\prime$ has only one geodesic boundary. 
This is a contradiction.
Hence, we have $N(X, \omega) \leq g$. 

Next, we assume that  $(X, \omega)\in \calhhyp(g-1, g-1)$.
If $N(X, \omega) \geq g+2$, 
then there exist  disjoint regular closed geodesics 
$\gamma_1, \dots, \gamma_{g+2}$ in $(X, \omega)$ that are not homotopic to each other. 
Since $(X, \omega)$ is a surface of genus $g$, 
$X-\left(\gamma_1 \cup \cdots \cup \gamma_{g+1}\right)$ has 
at least 
three connected components.
Since  $(X, \omega)$ has only two singular points, one of the connected components $X^\prime$ 
 is a translation surface that has no singular points. 
By the same argument as above, 
%$X^\prime$ is a Euclidean cylinder.
%Moreover,  the boundary of $X^\prime$ consists of some $\gamma_i$ and $\gamma_j$.
%Thus, $\gamma_i$ and $\gamma_j$ are homotopic to each other.
this is a contradiction.
Hence, we have $N(X, \omega) \leq g+1$. 
}

\subsection{The idea of lower estimates of $N(X, \omega)$}
Hereafter, we prove the inequality 
\begin{align}\label{lower_bound}
N(X, \omega) \geq 
\begin{cases}
g & \left((X, \omega)\in \calhhyp(2g-2)\right)  \\
g+1 & \left((X, \omega)\in \calhhyp(g-1, g-1) \right).
\end{cases} 
\end{align}
To do this, we consider the quotient of 
$(X, \omega) \in \calhhyp (2g-2)\cup \calhhyp(g-1, g-1)$ 
by its hyperelliptic involution $\tau$.
By definition, we have $\tau^\ast \omega=-\omega$.
Thus, the holomorphic quadratic differential $\omega^2$ 
induces a meromorphic quadratic differential $q$ on $Y=X/\left< \tau \right>$
via the natural projection $\varphi : X \to Y$. 
Then we have the following.

%If $p_0$ is a unique zero of $\omega$, the angle around $p_0$ 
%is $4g-2$.
%The image $\varphi(p_0)$ is a unique zero of $q$ and the order is $2g-3$
%The poles of $q$ are the imagees of Weierstrass points of $X$ 
%
%
%
%
%
%the hyperelliptic involution  $\tau$ is locally realized as 
%a $\pi$-rotation with respect to 
%the Euclidean structure of $(X, \omega)$. 
%This implies that the singular Euclidean structure on $(X, \omega)$ 
%induces a singular Euclidean structure on $\overline{Y}=X/\left< \tau \right>$.
%The singular Euclidean structure is 
%
%
%
%
%By the projection $\varphi : X \to \overline{Y}=X/\left< \tau \right>$,  
%the quadratic differential $\omega^2$ on $X$ induces 
%a holomorphic quadratic differential $q$ 
%
%
%
%$$ and 
% the  of $X$. 
%By the projection $\varphi : X \to \overline{Y}=X/\left< \tau \right>$,  
%the quadratic differential $\omega^2$ on $X$ induces 
%a holomorphic quadratic differential $q$ 
%on $Y=\overline{Y}-\varphi ({\rm Fix} (\tau))$.

\prop{
The Riemann surface $Y=X/\left< \tau \right>$ is a compact surface of genus $0$.
If $(X, \omega) \in \calhhyp (2g-2)$, 
the quadratic differential $q$ has a unique zero of order $2g-3$ 
and $2g+1$ simple poles.
If $(X, \omega) \in \calhhyp (g-1, g-1)$, 
the quadratic differential  $q$ has a unique zero of order $2g-2$ and 
$2g+2$ simple poles.
}

%
% Let $(X, \omega)$ be a translation surface of genus $g \geq 2$.
% We assume that $(X, \omega )$ is in $\calhhyp(2g-2)$ or $\calhhyp(g-1, g-1)$. 
%Let $\tau$ be the hyperelliptic involution of $(X, \omega)$.
%Then $Y=X/\left< \tau \right>$ is a Riemann

\pf{
Since $\tau$ is the hyperelliptic involution, $\tau$ has $2g+2$ fixed points, 
and $Y=X/\left< \tau \right>$ is a compact surface of genus $0$.
Suppose that $(X, \omega) \in \calhhyp (2g-2)$.
Then the unique zero of $\omega$ is a fixed point of  $\tau$.
As the angle around $p_0$ on $(X, \omega)$ is $(4g-2)\pi$, 
the angle  around $\varphi(p_0)$ on $(Y, q)$ is $(2g-1)\pi$.
Therefore, $\varphi(p_0)$ is a zero of $q$ of order $2g-3$.
The angle of every other fixed point $p$ of $\tau$ is $2\pi$.
Thus, the angle around $\varphi(p)$ on $(Y, q)$ is $\pi$.
This means that the point is a simple pole of $q$.
Every point of $X \setminus {\rm Fix} (\tau)$ 
is mapped to a point on $Y$ that is not a zero nor a pole of $q$.
Next, we assume that $(X, \omega) \in \calhhyp (g-1, g-1)$. 
Then, all fixed points of $\tau$ are mapped to simple poles of $q$.
Two zeros of $\omega$  are mapped to a point on $(Y, q)$ that is zero of $q$ of  order $2g-2$.
Every point of $X \setminus {\rm Fix} (\tau)$ 
is neither mapped to a zero nor a pole of $q$. 
%Now, we obtain the claim.
}

The inequality (\ref{lower_bound}) is proved by giving 
disjoint regular simple closed geodesics 
$\gamma_1, \dots, \gamma_g$ 
(resp. $\gamma_1, \dots, \gamma_{g+1}$) 
in $(X, \omega) \in \calhhyp(2g-2)$
(resp. $(X, \omega) \in \calhhyp(g-1, g-1)$)
each of which is invariant under the hyperelliptic involution $\tau$ of $(X, \omega)$. 
Let $\gamma$ be a regular closed geodesic in 
$(X, \omega)$ that is invariant under $\tau$. 
Since every orientation reversing isometry of $S^1$ has two fixed points, 
the image $\varphi(\gamma)$  is a saddle connection on $(Y, q)$ 
that connects distinct simple poles of $q$.
Conversely, if $s$ is such a saddle connection on $(Y, q)$, 
the preimage $\varphi^{-1}(s)$ is  a regular closed geodesic on $(X, \omega)$ 
that is invariant under $\tau$.
Therefore, we find  such saddle connections on $(Y, q)$ that are disjoint to each other.
%In our case, $(Y, q)$ is always a flat surface of genus $0$ such that $q$ has a unique zero
%and some simple poles of $q$. 
To do this, we give a triangulation of $(Y, q)$ by ``restricted triangles'' defined by Definition \ref{restricted_triangle}.
All vertices of this triangulation are the zero of $q$ and 
all poles of $q$ are midpoints of some edges.
%This triangulation easily gives us saddle connections that we want. 

%To do this, we give a triangulationf of $(Y, q)$.
%Howver, we can not always obtain such saddle connections from    

\subsection{Geometry of flat surfaces of genus $0$}
In this subsection, we give some properties of flat surfaces of genus $0$. 
Let $(Y, q)$ be a flat surface of genus $0$.
Note that $(Y, q)$ must have some simple poles.
%Then $q$ has some simple poles.
Let $s^\prime : [0,1] \to (Y, q)$ be a geodesic segment connecting a singular point to a simple pole of $q$.
Then 
$$s(t)=
\begin{cases}
s^\prime (2t)  & (0 \leq t \leq \frac{1}{2}) \\
s^\prime(2-2t) & (\frac{1}{2} \leq t \leq 1)
\end{cases}
$$
is a closed curve on $(Y, q)$ 
such that its image is a segment.

\defi{[Returning geodesic]
We call a closed curve $s$ as above a \emph{returning geodesic}.
That is, a returning geodesic $s: [0, 1] \to (Y, q)$ is a closed curve from a singular point 
of $(Y, q)$ to itself such that  $s\left(\frac{1}{2}\right)$ is a simple pole of $q$, 
$s|_{\left[0, \frac{1}{2}\right]}$ is a saddle connection, 
and $s(t)=s(1-t)$ for all $t\in [0,1]$.
}

The next lemma is important for our proof.

\lem{\label{return}
Let $(Y, q)$ be a flat surface of genus $0$.
Let $s: [0,1] \to (Y, q)$ be a curve satisfying the following:
\begin{itemize}
\item the absolute value of the derivative $|s^\prime|$ is constant, 
\item $s(0)$ is a singular point and $s\left(\frac{1}{2} \right)$ is a simple pole, 
\item $s|_{\left[0, \frac{1}{2} \right]}$ is a local isometry 
such that  $s\left(\left(0, \frac{1}{2} \right)\right)$ contains no singular points, 
\item $s|_{\left[\frac{1}{2}, 1\right]}$ is a local isometry and  
\item there exists  $\varepsilon>0$ 
such that
$s \left(1-t\right)=s \left(t\right)$ 
for all $t\in \left(\frac{1}{2}-\varepsilon, \frac{1}{2}+\varepsilon \right)$.
%
%\item the absolute value of the derivative $|s^\prime|$ is constant, 
%\item $s(0)$ is a singular point and $s\left(\frac{1}{2} \right)$ is a simple pole, 
%\item $s\left(\left[\frac{1}{2}, 1 \right]\right)$ is a geodesic and  
%\item there exists a chart $(U, z)$ of $(Y, q)$ around $s\left(\frac{1}{2} \right)$ 
%and a sufficiently small $\varepsilon>0$
%such that $z\circ s |_{ \left( \frac{1}{2}-\varepsilon, \frac{1}{2}+\varepsilon\right)} : \left( \frac{1}{2}-\varepsilon, \frac{1}{2}+\varepsilon\right) \to \mathbb{C}$ is a segment
\end{itemize}
% and 
%. 
%Assume that 
%
%
%$s(0)$ is a singular point of $(Y, q)$.
%If  $s\left(\left(0, \frac{1}{2} \right)\right)$ contains no singular points 
%and $s\left(\frac{1}{2} \right)$ is a simple pole of $(Y, q)$, 
Then $s(t)=s(1-t)$ for all $t\in [0, 1]$. 
Especially, $s\left(\left(\frac{1}{2}, 1 \right)\right)$ contains no singular points
and $s(1)$ coincides with the singular point $s(0)$.
}

\pf{
Since  
$s|_{\left[0, \frac{1}{2} \right]}$
and
$s|_{\left[\frac{1}{2}, 1 \right]}$
are local isometries and 
$s\left(\frac{1}{2} \right)$ is a simple pole, $s$ is a returning geodesic. 
Thus, we obtain the claim.
}

\prop{\label{embedding}
Let $(Y, q)$ be a flat surface of genus $0$.
Let $s : [0,1] \to (Y, q)$ be a returning geodesic.
Given a segment $s_0$ in $\mathbb{C}$ such that 
$s_0$ is parallel to $s$ and  the length $|s_0|$ equals the length of $s$. 
Choose 
 a non-zero vector $\vv$ that is not parallel to $s_0$. 
Then there exists an  orientation preserving  local isometric embedding 
$\rho_0 : D \to (Y, q)$ satisfying the following:
\begin{enumerate}
\item $D$ is an open parallelogram 
such that 
one of the sides is $s_0$ and another side is parallel to $\vv$,
%\item $\rho$ is a local isometric embedding,
\item $\rho_0(D)$ contains no singular points,
\item $\rho_0$ can be extended to an immersion $\rho_0 : \overline{D} \to (Y, q)$, 
\item $s_0$ is mapped to $s$ via $\rho_0$,
\item the other side that is parallel to $s_0$ is mapped to an open segment containing a singular point.
%\item the sides that are parallel to $\vv$ are mapped to  
%segments containing singular points
% in their interiors.
\end{enumerate}
}

\pf{
For sufficiently small $\varepsilon >0$, 
there exists  an immersion $\rho_\varepsilon : D_\varepsilon \to (Y, q)$ 
satisfying the same conditions as $(1), (2), (3), (4)$ 
and the distance between $s_0$ and its opposite side 
%two sides that are parallel to $\vv_s$ 
is $\varepsilon$.
%
%For sufficiently small $\varepsilon >0$, the $\varepsilon$-neighborhood of $s$ 
%contains no singular points.
%Thus there exists  an immersion $\rho_\varepsilon : D_\varepsilon \to (Y, q)$ 
%satisfying the same conditions as $(1), (2), (3), (4)$.
%%Then we have $D_\varepsilon \subset D_{\varepsilon^\prime}$ if $\varepsilon < \varepsilon^\prime$.
%
We show that there exists the maximum of $\varepsilon>0$ such that the above 
$\rho_\varepsilon$ exists.
If not, there exists 
a local isometric immersion  from 
an open band with infinite area into $(Y, q)$. 
Since the area of $(Y, q)$ is finite, 
the image of this embedding must be a cylinder.
Since the angle around $s \left(\frac{1}{2} \right)$ is $\pi$, this is a contradiction.
Let $\varepsilon_0>0$ be the maximal of $\varepsilon>0$ such that the above 
$\rho_\varepsilon$ exists.
Then $\rho_0=\rho_{\varepsilon_0}$ also satisfies the condition $(5)$. 
%
%The maximal $\varepsilon_0$ exists since if not, there exists 
%an isometrically embedding from 
%an open band with infinite area into $(Y, q)$. 
%However the area of $(Y, q)$ is finite.
}

%\lem{\label{recurrence}
%Let $(Y, q)$ be a flat surface of genus $0$ such that $q$ has a unique zero.
%Choose a returning geodesic  $s$ on $(Y, q)$. 
%
%
%Let $\vv$ be a unit vector that is not parallel to $s$ and 
%$\left\{ F_t : t \geq 0 \right\}$ 
%the flow of direction $\vv$.
%Then there exists $t_0>0$ such that $F_{t_0}(s)$ contains a singular point of $q$.
%}
%
%
%\pf{
%If not,  there exists an isometrically embedding from 
%a band with infinite area into $(Y, q)$.
%However, the area of $(Y, q)$ is finite.
%This is a contradiction. 
%}

%\remark{
%Lemma \ref{recurrence} is also true if the segment $s$ contains 
%a  zero of order $-1$ in its interior.
%}

\defi{[$(s_0, \vv)$-restricted triangle and left (right) strongly $(s_0, \vv)$-restricted triangle]
\label{restricted_triangle}
 Let $\Delta$ be a Euclidean triangle in $\mathbb{C}$
with a side $s_0$. 
Let $\vv$ be a unit vector that is not parallel to $s_0$.
Assume that $P$ is the vertex of $\Delta$ that is opposite to $s_0$.
\begin{enumerate}
\item 
%\begin{enumerate}
%\item 
The triangle $\Delta$ is a $(s_0, \vv)$-restricted triangle
if the Euclidean line $l$ that  passes through $P$ 
and is of direction $\vv$ 
intersects with $s_0$. 
\item 
A left (resp. right) side of a $(s_0, \vv)$-restricted triangle $\Delta$ 
is the side of $\Delta$ that is at the left (resp. right) when  
we put $\Delta$ so that $s_0$ is horizontal and $P$ 
is above $s_0$.

\item A left (resp. right) strongly $(s_0, \vv)$-restricted triangle $\Delta$ 
is a $(s_0, \vv)$-restricted triangle 
whose left (resp. right)  side is not parallel to $\vv$. 
%
%\item A strongly $(s_0, \vv)$-restricted triangle $\Delta$ 
%is a $(s_0, \vv)$-restricted triangle  such that both 
%left side and right side are not parallel to $\vv$.
\end{enumerate}

%If the line $l$ intersect with the interior of $s_0$, 
%the triangle $\Delta$ is called a strongly  $(s_0, \vv)$-restricted triangle.
%\item 
%
%Let $\Delta$ be a $(s, \vv)$-restricted  triangle.
%Let $\rho : \Delta \to (Y, q)$ be an immersion into a flat surface $(Y, q)$ 
%such that $\rho|_{{\rm Int}(\Delta)}$ is local isometric embedding. 
%Then 
%\end{enumerate}
}

We prove the following theorem.

\thm{\label{triangle}
Let $(Y, q)$ be a flat surface of genus $0$
 such that $q$ has a unique zero $p_0$.
Let  $s$ be a returning geodesic on $(Y, q)$ that starts from the point $p_0$.
Let $\vv$ be a unit vector that is not parallel to $s$. 
Then there exists  
a (closed) triangle $\Delta$ in $\mathbb{C}$ 
with sides $s_0, s_1, s_2$
and an  orientation preserving  immersion 
$\rho : \Delta \to (Y, q)$ 
 satisfying the following: 
\begin{enumerate}
\item $s_0$ is parallel to $s$ and $\Delta$ is a left (resp. right) strongly $(s_0, \vv)$-restricted triangle, 
\item $\rho|_{{\rm Int}(\Delta)}$ is a local isometric embedding,
\item $\rho({\rm Int}(\Delta))$ contains no singular points, 
\item every vertex of $\Delta$ is mapped to $p_0$,
\item $s=\rho \circ s_0$ and 
\item $\rho \circ s_i$ is a saddle connection or a returning geodesic for each $i=1, 2$.
%\item If $s_i$ is parallel to $\vv$ for some $i=1,2$, then $\rho \circ s_i$ is a returning geodesic.
\end{enumerate}
}

We prove Theorem \ref{triangle} 
in the case where $\Delta$ is a left strongly $(s_0, \vv)$-restricted triangle. 
The right strongly $(s_0, \vv)$-restricted triangle case is proved in the same way. 
Theorem \ref{triangle} is proved by the following five Lemmas. 

\begin{figure}[ht!]
\labellist
%\small
\hair 0pt
\pinlabel $O$  at 210 -15
\pinlabel $s_0$  at 300 -10
\pinlabel $a$  at 255 165
\pinlabel $b$  at 40 105
\pinlabel $c$  at 500 105
\pinlabel $\vv$  at 10 50
\pinlabel $P_1$  at -20 5
\pinlabel $P_2$  at 435 5
\pinlabel $P_3$  at 530 145
\pinlabel $P_4$  at 50 145
 \endlabellist
\centering
 \includegraphics[scale=0.4]{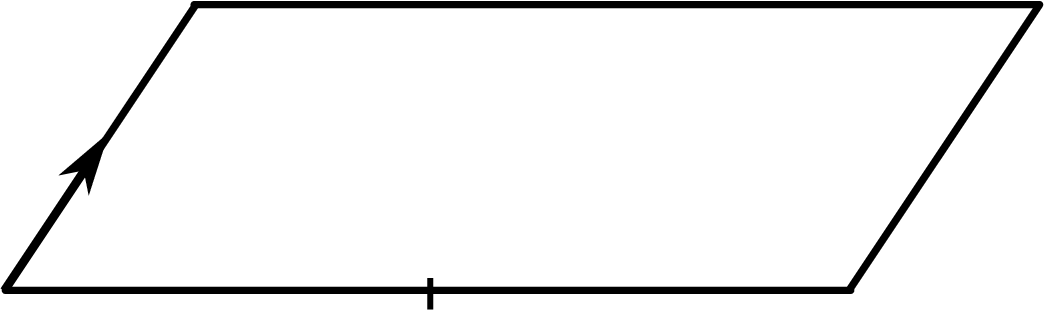} 
\caption{The parallelogram $D$.}
\label{immersionD}
\end{figure}

%%%%%%%
Let $(Y, q)$ be a flat surface of genus $0$
 such that $q$ has a unique zero $p_0$.
Let  $s$ be a returning geodesic on $(Y, q)$ that starts from the point $p_0$.
Given $s_0$, $\vv$, $D$, and  $\rho_0: D \to (Y, q)$ as in Proposition \ref{embedding}.
%Let $s_0$ be a side of $D$ satisfying 
Then we have $\rho_0 \circ s_0=s$. 
We label all other sides of $D$ by $a, b, c$ 
so that the counterclockwise order of the sides is $s_0, c, a, b$. 
%and $a$ the opposite side of $s_0$. 
%We denote by $P_1$ the common point of $s_0$ and $b$.
%We also  denote by $P_2$ the common point of $s_0$ and $c$ (See 
We label the vertices of $D$ by $P_1, P_2, P_3, P_4$ as in 
Figure \ref{immersionD}. 
We also give the label $O$ to the midpoint of $s_0$.
%).

We define the following terminologies.
\begin{itemize}
\item For a point $P\in \mathbb{C}$, we denote by $r_P$ 
 the point reflection in $P$.
\item A point $P \in a$ is a left point of $a$ if      $P_4P <\frac{1}{2}|a|$. 
\item A point $P \in a$ is a right point of $a$ if    $P_4P >\frac{1}{2}|a|$. 
\item A point $P \in b$ is a lower point of $b$ if   $P_1P <\frac{1}{2}|b|$.
\item A point $P \in b$ is an upper point of $b$ if $P_1P >\frac{1}{2}|b|$.
\item A point $P \in c$ is a lower point of $c$ if   $P_2P <\frac{1}{2}|c|$. 
\item A point $P \in c$ is an upper point of $c$ if $P_2P >\frac{1}{2}|c|$. 
\end{itemize}

\lem{\label{trivial}
If the interior of the side $a$ contains a point $P$ that is mapped via $\rho_0$ 
to $p_0$, there exists a triangle $\Delta$ in $\mathbb{C}$ 
with sides $s_0, s_1, s_2$
and an  orientation preserving  immersion 
$\rho : \Delta \to (Y, q)$ 
 satisfying all conditions as in Theorem \ref{triangle}. 
}

\pf{
Setting $\Delta=\triangle P_1P_2P$ and $\rho=\rho_0|_\Delta$, we obtain the claim.
}

\lem{\label{c_sing}
Assume that the interior of the side $a$ does not  contain a point that is mapped via $\rho_0$ 
to $p_0$ and 
$\rho_0(c-\left\{P_2\right\})$ contains a singular point.
Then, there exist a triangle $\Delta$ in $\mathbb{C}$ 
with sides $s_0, s_1, s_2$
and an orientation preserving  immersion 
$\rho : \Delta \to (Y, q)$ 
 satisfying all conditions as in Theorem \ref{triangle}. 
}

\begin{figure}[ht!]
\labellist
%\small
\hair 0pt
\pinlabel $O$  at 210 -15
\pinlabel $s_0$  at 300 -10
\pinlabel $a$  at 255 165
\pinlabel $b$  at 15 70
\pinlabel $c$  at 475 70
%\pinlabel $\vv$  at 10 50
\pinlabel $P_1$  at -20 5
\pinlabel $P_2$  at 435 5
\pinlabel $P_3$  at 530 135
\pinlabel $P_4$  at 50 135
\pinlabel $P$  at 500 105
\pinlabel $r_P(P_2)$  at 590 210
\pinlabel $R$  at 200 165
\pinlabel $c^\prime$  at 545 180
 \endlabellist
\centering
 \includegraphics[scale=0.4]{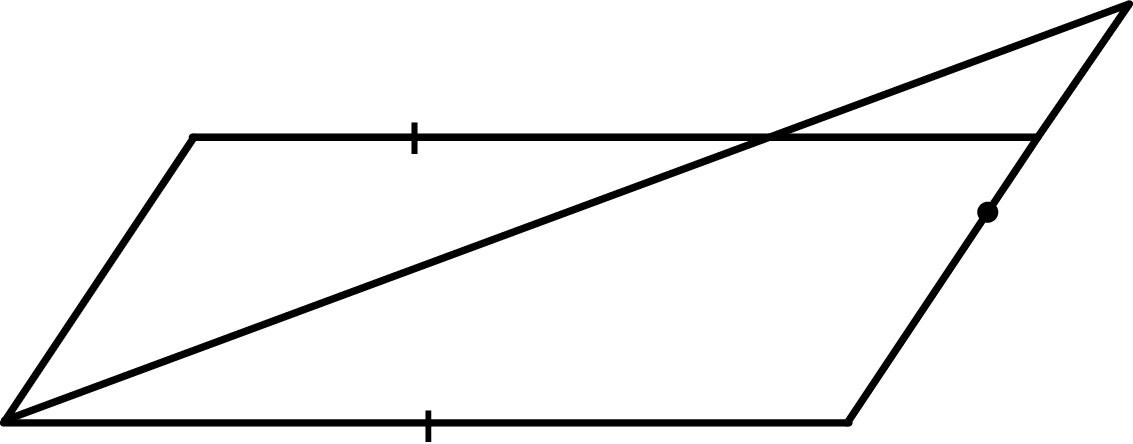} 
\caption{Case (ii) with the assumption that  $R$ is the left point of $a$}
\label{case_two_left}
\end{figure}

\pf{
Let $P$ be the point in $c-\left\{P_2\right\}$ such that 
$P$ is closest to $P_2$ in 
the points  of  $c-\left\{P_2\right\}$  
that are mapped to singular points of $(Y, q)$. 
If $\rho_0(P)=p_0$, then 
$\Delta=\triangle P_1P_2P$ and $\rho=\rho_0|_\Delta$ 
satisfies all conditions as in Theorem \ref{triangle}. 
Hereafter, we assume that $\rho_0(P)$ is a pole of $q$. 
%Let $r_{P}$ be the point reflection in $P$.

Case (i) Suppose that $P$ is a lower point or the midpoint of $c$. 
By Lemma \ref{return}, $Q=r_{P}(P_2) \in c$ is 
a point whose image via $\rho_0$ is $p_0$ (see Figure \ref{case_one}).
Moreover, there exist no points whose images are singular points between $P$ and $Q$. 
Setting $\Delta=\triangle P_1P_2Q$ and $\rho=\rho_0|_\Delta$, 
$\Delta$ and $\rho$  
satisfy the conditions 
as in Theorem \ref{triangle}.

\begin{figure}[ht!]
\labellist
%\small
\hair 0pt
\pinlabel $O$  at 210 -15
\pinlabel $s_0$  at 300 -10
\pinlabel $a$  at 255 165
\pinlabel $b$  at 40 105
\pinlabel $c$  at 500 105
%\pinlabel $\vv$  at 10 50
\pinlabel $P_1$  at -20 5
\pinlabel $P_2$  at 435 5
\pinlabel $P_3$  at 530 145
\pinlabel $P_4$  at 50 145
\pinlabel $P$  at 450 45
\pinlabel $Q$  at 475 80
 \endlabellist
\centering
 \includegraphics[scale=0.4]{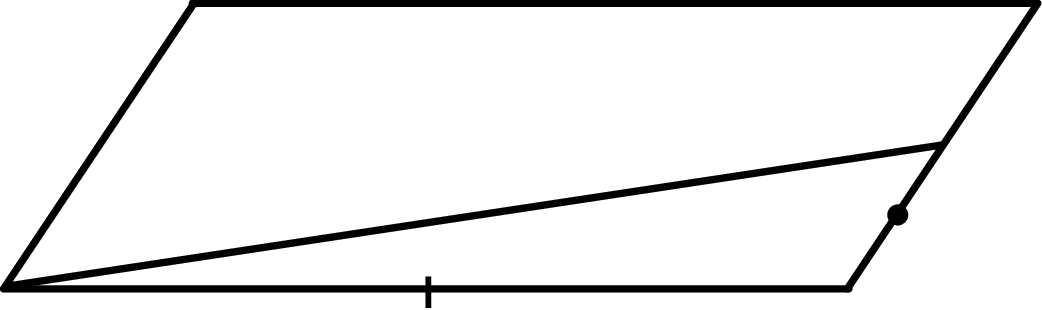} 
\caption{Case (i).}
\label{case_one}
\end{figure}

%%%%%%%%%%%figure_a
Case (ii)  Suppose that $P$ is an upper point of $c$.
By assumption, there exists  a point $R \in \interior{a}$ 
whose image is a simple pole of $q$. 
Then, $R$ is not a right point of $a$. 
If $R$ is a right point of $a$, by Lemma \ref{return}, 
$r_R \circ r_P(P_2)$ is 
an interior point of $D$.
However, $r_R \circ r_P(P_2)$ is mapped to $p_0$. 
This is a contradiction.
If  $R$ is the midpoint of $a$, 
we set $\Delta=\triangle P_1P_2 r_P(P_2)$ and 
$$\rho(z)=
\begin{cases}
\rho_0 (z)   & (z \in \Delta \cap D) \\
\rho_0 \circ r_R(z) & (z\in \Delta-D)
\end{cases}
$$
(see Figure \ref{case_two_mid}).
Then, $\Delta$ and $\rho$ satisfy the conditions 
as in Theorem \ref{triangle}.

\begin{figure}[ht!]
\labellist
%\small
\hair 0pt
\pinlabel $O$  at 210 -15
\pinlabel $s_0$  at 300 -10
\pinlabel $a$  at 255 165
\pinlabel $b$  at 15 70
\pinlabel $c$  at 475 70
%\pinlabel $\vv$  at 10 50
\pinlabel $P_1$  at -20 5
\pinlabel $P_2$  at 435 5
%%\pinlabel $P_3$  at 530 145
%%\pinlabel $P_4$  at 50 145
\pinlabel $P$  at 500 105
\pinlabel $r_P(P_2)$  at 590 210
\pinlabel $R$  at 300 165
 \endlabellist
\centering
 \includegraphics[scale=0.4]{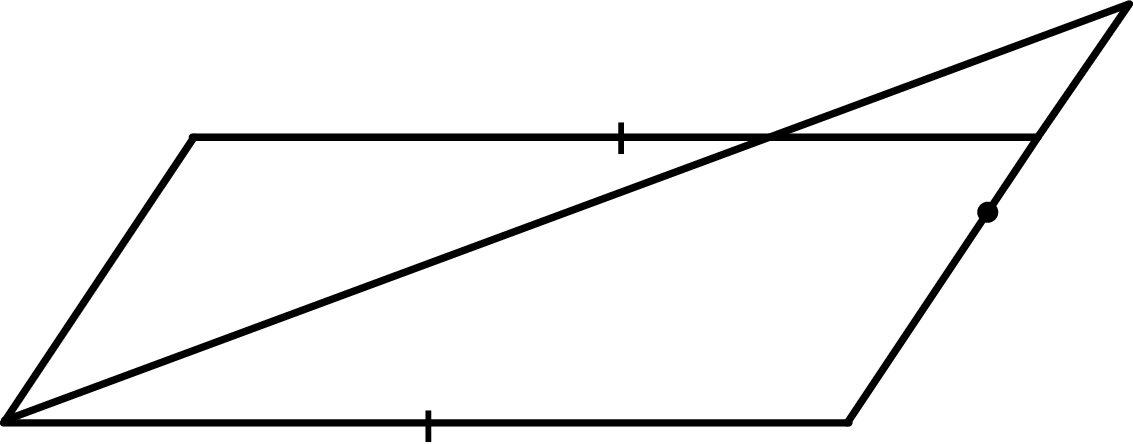} 
\caption{Case (ii) with the assumption that  $R$ is the midpoint of $a$.}
\label{case_two_mid}
\end{figure}

%%%%%%%%%%%figure_b
Finally, suppose that  $R$ is a left point of $a$.
Let $E=\triangle P_1P_2r_P(P_2) - D$ and 
$c^\prime$ the side of $E$ that is parallel to $c$ 
(see Figure \ref{case_two_left}). 
If $\rho_0$ can be extended to a local isometry 
from $\interior{D\cup E}$ to $(Y, q)$,  then 
$\Delta=\triangle P_1P_2r_P(P_2) $ and 
 $\rho=\rho_0 |_\Delta$ 
 satisfy the conditions as in Theorem \ref{triangle}.  
We now assume that $\rho_0$ cannot be extended to $\interior{D\cup E}$ 
as an isometry.
Let $\theta_0=\angle P_3 P_1 P_2$ 
and $\theta_1=\angle r_P(P_2) P_1 P_2$. 
For each $\theta$ ($\theta_0 <\theta <\theta_1$), 
let $l_\theta$ be the segment  starting from $P_1$ to a point of $c^\prime$ such that 
the angle between $P_1P_2$ and $l_\theta$ is $\theta$.
By assumption, 
there exists the maximal $\Theta$
such that 
$\rho_0$ can be extended to a local isometry 
from $ {\rm Int} 
\left( 
D\cup  \bigcup_{\theta_0 <\theta<\Theta}  l_\theta
 \right)$ 
to $(Y, q)$.
Then, $l_\Theta$ contains points that are mapped to singular points.
Let $S$ be the point that is closest to $s_0$ in the points.
We show that $S$ is mapped to the zero $p_0$ of $q$.
Suppose that  $S$ is mapped to a pole of $q$.
Let $T_1$ be the intersection of segments $P_1 r_s(P_1)$ and $P_2r_P(P_2)$. 
Let $T_2$ be the point in $P_1 r_S(P_1)$ 
such that $r_P(T_2)$ is the intersection of $P_1 P_2$ and $r_P\left(T_1 r_S(P_1) \right)$. 
We define a curve $s: P_1r_S(P_1) \to (Y, q)$ by 
 $$s(z)=
\begin{cases}
\rho_0 (z)   & \left(z \in P_1 T_1 \right) \\
\rho_0 \circ r_P(z) & \left(z\in T_1 T_2 \right)\\
\rho_0 \circ r_O\circ r_P(z) & \left(z\in T_2 r_S(P_1) \right).
\end{cases}
$$
By Lemma \ref{return},  
we have $\rho_0 \circ r_O \circ  r_P \circ r_S(P_1)=p_0$. 
However, we have $ r_O \circ  r_P\circ r_S(P_1) \in D$.
This is a contradiction. 
Therefore, $S$ is mapped to $p_0$. 
Setting $\Delta=\triangle P_1 P_2 S$ and $\rho=\rho_0|\Delta$, 
then $\Delta$ and $\rho$ satisfy the conditions 
as in Theorem \ref{triangle}.
}

%%%%%%%%%%%%%%%%%%c

\lem{\label{mid}
Assume that the interior of the side $a$ does not  contain a point that is mapped 
to $p_0$, 
the midpoint $R$ of $a$ is mapped to a pole of $q$ 
and 
$\rho_0(c-\left\{P_2\right\})$ contains no singular points.
%and $\rho_0(b-\left\{P_1\right\})$ contains a singular point
Then, there exist a triangle $\Delta$ in $\mathbb{C}$ 
with sides $s_0, s_1, s_2$
and an orientation preserving  immersion 
$\rho : \Delta \to (Y, q)$ 
 satisfying all conditions as in Theorem \ref{triangle}. 
}

\begin{figure}[ht!]
\labellist
%\small
\hair 0pt
\pinlabel $O$  at 210 -15
\pinlabel $s_0$  at 300 -10
\pinlabel $a$  at 255 165
\pinlabel $b$  at 15 70
\pinlabel $c$  at 475 70
%\pinlabel $\vv$  at 10 50
\pinlabel $P_1$  at -20 5
\pinlabel $P_2$  at 435 5
%\pinlabel $P_3$  at 530 145
%\pinlabel $P_4$  at 50 145
\pinlabel $P$  at 50 105
\pinlabel $r_R(P)$  at 550 180
\pinlabel $R$  at 300 165
 \endlabellist
\centering
 \includegraphics[scale=0.4]{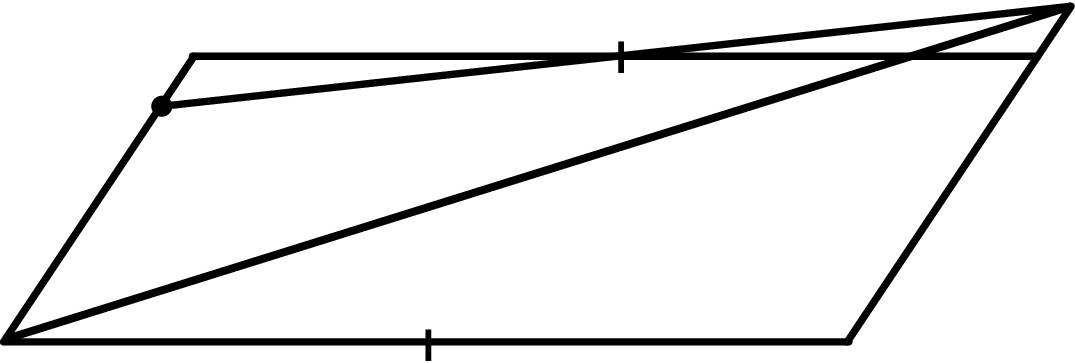} 
\caption{The case of Lemma \ref{mid}.}
\label{b_sing_R_mid}
\end{figure}

\pf{
Let $P$ be the point in $b$ such that 
$P$ is farthest from $P_1$ in 
the points  of  $b$  
that are mapped to singular points of $(Y, q)$. 
Note that $P$ may possibly coincide with $P_1$.
We first prove  that  $P$ is mapped to $p_0$. 
Suppose that $P$ is mapped to a pole of $q$.
If $P$ is a lower point or the midpoint of $b$, 
by Lemma \ref{return}, 
$r_P(P_1)\in b$ is mapped to a singular point.
This contradicts the definition of $P$. 
Assume that $P$ is an upper point of $b$.
We show that $P$ is mapped to $p_0$. 
Suppose that $P$ is mapped to a pole of $q$.
Let $P^\prime$ be the point 
that is closest to $P$  in the points of  $P_1 P$ corresponding to singular points 
and $P^\prime \neq P$.
Let $s$ be a curve $s: P^\prime r_P(P^\prime) \to (Y, q)$ defined by 
 $$s(z)=
\begin{cases}
\rho_0 (z) & \left(z\in  P^\prime P_4  \right)\\
\rho_0 \circ r_R (z)   & \left(z \in  P_4 r_P(P^\prime) \right).
\end{cases}
$$
By Lemma \ref{return}, we have  $s(r_P(P^\prime) )=s(P^\prime)$. 
Hence, $s(r_P(P^\prime) )=\rho_0 \circ r_R  (r_P(P^\prime) ) =\rho_0 (r_R r_P (P^\prime)) $ is a singular point.
However, 
$r_R r_P (P^\prime)$ is in $c-\left\{P_2\right\}$.
This is a contradiction.
Therefore, $P$ is mapped to $p_0$. 
%Since $\rho_0( P^\prime  \right\})$ contains no singular points, 
%$s(r_P(P_1) P_4)$ contains no singular points. 
%Especially, $s(r_P(P_1))$ is not a singular point.
%By the definition of $P$, 
%$s(P_4 P)$ contains no singular points. 
%By Lemma \ref{return},  we have 
%$s(r_P(P_1))=s(P_1)=\rho_0(P_1)=p_0$.  
%%However, 
%%$r_R r_P(P_1)$ is in $c-\left\{P_2\right\}$.
%%However,  
%%since 
%%$r_R r_P(P_1)$ is identified with $P_1$ via 
%%by 
%%mapped to $p_0$.
%%by the same argument as in the last part of 
%%the proof of Lemma \ref{c_sing}, 
%%%Lemma \ref{return}, 
%%$r_R r_P(P_1)$ is mapped to $p_0$.
%%However, 
%This is a contradiction. 
%Thus, $P$ is mapped to $p_0$. 
Set $\Delta=\triangle P_1P_2r_R(P)$ and 
$$\rho(z)=
\begin{cases}
\rho_0 (z)   & (z \in \Delta \cap D) \\
\rho_0 \circ r_R(z) & (z\in \Delta-D)
\end{cases}
$$
(see Figure \ref{b_sing_R_mid}).
Then, $\Delta$ and $\rho$ satisfy the conditions.
}

\lem{\label{left}
Assume that the interior of the side $a$ does not  contain a point that is mapped 
to $p_0$, 
the left point $R$ of $a$ is mapped to a pole of $q$,  
and 
$\rho_0(c-\left\{P_2\right\})$ contains no singular points. 
%and $\rho_0(b-\left\{P_1\right\})$ contains a singular point 
Then, there exist a triangle $\Delta$ in $\mathbb{C}$ 
with sides $s_0, s_1, s_2$
and an orientation preserving  immersion 
$\rho : \Delta \to (Y, q)$ 
 satisfying all conditions as in Theorem \ref{triangle}. 
}

\begin{figure}[ht!]
\labellist
%\small
\hair 0pt
\pinlabel $O$  at 210 -15
\pinlabel $s_0$  at 300 -10
\pinlabel $a$  at 400 165
\pinlabel $b$  at 15 70
\pinlabel $c$  at 475 70
%\pinlabel $\vv$  at 10 50
\pinlabel $P_1$  at -20 5
\pinlabel $P_2$  at 435 5
\pinlabel $P_3$  at 530 150
\pinlabel $P_4$  at 60 150
\pinlabel $P$  at 35 105
\pinlabel $r_R(P)$  at 265 200
\pinlabel $R$  at 150 170
\pinlabel $b^\prime$  at 170 115
 \endlabellist
\centering
 \includegraphics[scale=0.4]{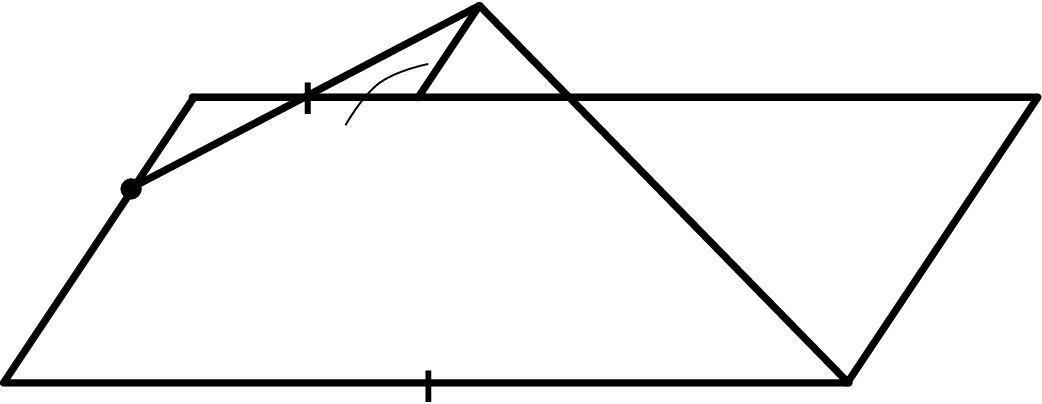} 
\caption{The case of Lemma \ref{left}.}
\label{b_sing_upper_R_left}
\end{figure}

\pf{
Let $P$ be the point in $b$ such that 
$P$ is farthest from $P_1$ in 
the points  of  $b$  
that are mapped to singular points of $(Y, q)$. 
Then, $P$ is mapped to $p_0$ by the same argument as in the proof of Lemma \ref{mid}. 
%Now we assume that $P$ is mapped to a pole of $(Y, q)$.
We make the same argument as the last argument of the proof of Lemma \ref{c_sing}.
Let $E=\triangle P_1P_2r_R(P) - D$ and 
$b^\prime = r_R(PP_4)$ 
(see Figure \ref{b_sing_upper_R_left}). 
If $\rho_0$ can be extended to a local  isometry 
from $\interior{D\cup E}$ to $(Y, q)$, 
then $\Delta=\triangle P_1P_2r_R(P) $ and 
 $\rho=\rho_0 |_\Delta$ 
 satisfy the conditions as in Theorem \ref{triangle}.  
We now assume that $\rho_0$ cannot be extended to $\interior{D\cup E}$ 
as a local  isometry.
Let $\theta_0=\angle P_1 P_2 r_R(P_4)$ 
and $\theta_1=\angle P_1 P_2 r_R(P)$. 
For each $\theta$ ($\theta_0 <\theta <\theta_1$), 
let $l_\theta$ be the segment  starting from $P_2$ to a point of $b^\prime$ such that 
the angle between $P_1P_2$ and $l_\theta$ is $\theta$.
By assumption, 
there exists the maximal $\Theta$
such that 
$\rho_0$ can be extended to a local  isometry 
from $D^\prime= {\rm Int} 
\left( 
D\cup  \bigcup_{\theta_0 <\theta<\Theta}  l_\theta
 \right)$ 
to $(Y, q)$.
Then, $l_\Theta$ contains points that are mapped to singular points except for $P_2$.
Let $S$ be the point that is closest to $s_0$ in the points.
If $S$ is mapped to a pole of $q$, 
%by Lemma \ref{return},  
$r_O \circ r_R \circ r_S(P_2)$ is mapped to $p_0$. 
However, we have $r_O \circ r_R\circ r_S(P_2)\in D$. 
This is a contradiction. 
Therefore, $S$ is mapped to $p_0$. 
Setting $\Delta=\triangle P_1 P_2 S$ and 
$$\rho(z)=
\begin{cases}
\rho_0 (z)   & (z \in \Delta \cap \overline{D^\prime}) \\
\rho_0 \circ r_R(z) & (z\in \Delta -\overline{D^\prime} ), 
\end{cases}
$$
%$\rho=\rho_0|_\Delta$,
$\Delta$ and $\rho$ satisfy the conditions 
as in Theorem \ref{triangle}.
}

\lem{\label{right}
Assume that the interior of the side $a$ does not  contain a point that is mapped 
to $p_0$,  
a right point $R$ of $a$ is mapped to a pole of $q$, and 
$\rho_0(c-\left\{P_2\right\})$ contains no singular points.
%and $\rho_0(b-\left\{P_1\right\})$ contains a singular point
%We also assume that 
%a right point $R$ of $a$ is mapped to a pole of $q$.
Then, there exist a triangle $\Delta$ in $\mathbb{C}$ 
with sides $s_0, s_1, s_2$
and an orientation preserving  immersion 
$\rho : \Delta \to (Y, q)$ 
 satisfying all conditions as in Theorem \ref{triangle}. 
}

\begin{figure}[ht!]
\labellist
%\small
\hair 0pt
\pinlabel $O$  at 210 -15
\pinlabel $s_0$  at 300 -10
%\pinlabel $a$  at 140 165
\pinlabel $b$  at 15 70
\pinlabel $c$  at 475 70
%\pinlabel $\vv$  at 10 50
\pinlabel $P_1$  at -20 5
\pinlabel $P_2$  at 435 5
\pinlabel $P_3$  at 530 150
\pinlabel $P_4$  at 60 150
%\pinlabel $P$  at 35 105
\pinlabel $r_R(P_2)$  at 420 280
\pinlabel $r_R(P_3)$  at 300 125
\pinlabel $R$  at 400 165
\pinlabel $c^\prime$  at 330 200
 \endlabellist
\centering
 \includegraphics[scale=0.4]{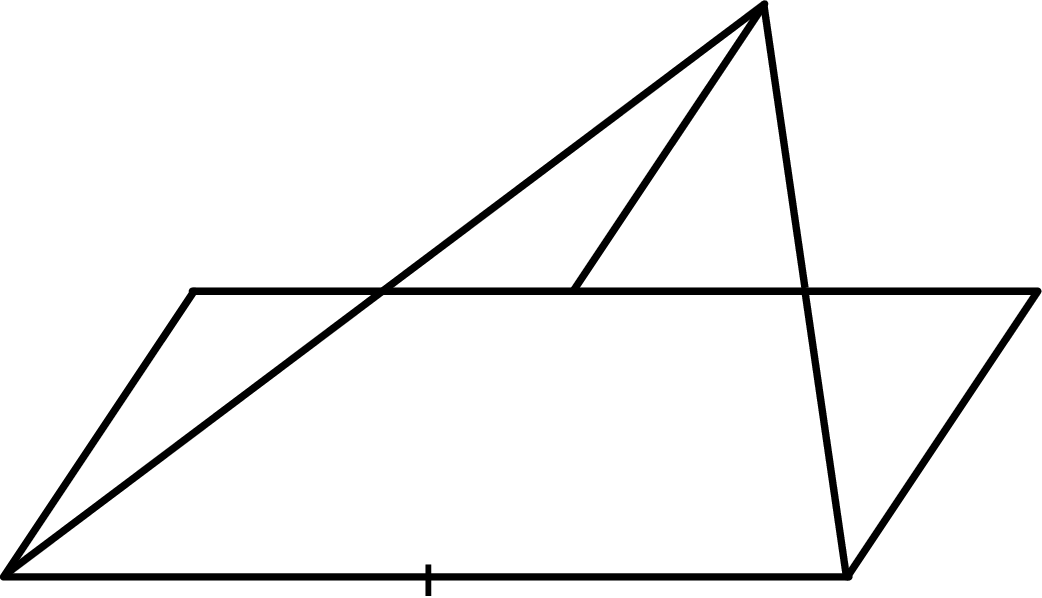} 
\caption{The case of Lemma \ref{right}.}
\label{R_right}
\end{figure}

\pf{
The proof is the same as the last argument of the proof of Lemma \ref{c_sing}.
Let $E=\triangle P_1P_2r_R(P_2) - D$ and 
$c^\prime = r_R(c)$ 
(see Figure \ref{R_right}). 
If $\rho_0$ can be extended to a local  isometry 
from $\interior{D\cup E}$ to $(Y, q)$, 
$\Delta=\triangle P_1P_2r_R(P_2) $ and 
 $\rho=\rho_0 |_\Delta$ 
 satisfy the conditions as in Theorem \ref{triangle}.  
We now assume that $\rho_0$ cannot be extended to $\interior{D\cup E}$ 
as a local isometry.
Let $\theta_0=\angle P_2 P_1 r_R(P_3)$ 
and $\theta_1=\angle P_2 P_1 r_R(P_2)$. 
For each $\theta$ ($\theta_0 <\theta <\theta_1$), 
let $l_\theta$ be the segment  starting from $P_2$ to a point of $c^\prime$ such that 
the angle between $P_1P_2$ and $l_\theta$ is $\theta$.
By assumption, 
there exists the maximal $\Theta$
such that 
$\rho_0$ can be extended to a local  isometry 
from $
D^{\prime\prime}= {\rm Int} 
\left( 
D\cup  \bigcup_{\theta_0 <\theta<\Theta}  l_\theta
 \right)$ 
to $(Y, q)$.
Then, $l_\Theta$ contains points that are mapped to singular points.
Let $S$ be the point that is closest to $s_0$ in the points.
If $S$ is mapped to a pole of $q$, 
%by Lemma \ref{return}, 
then  $r_O \circ r_R \circ r_S(P_1)$ is mapped to $p_0$. 
However, we have $r_O \circ r_R \circ r_S(P_1)\in D$. 
This is a contradiction. 
Therefore, $S$ is mapped to $p_0$. 
Setting $\Delta=\triangle P_1 P_2 S$ and 
%$\rho=\rho_0|_\Delta$, 
$$\rho(z)=
\begin{cases}
\rho_0 (z)   & (z \in \Delta \cap \overline{D^{\prime\prime}}) \\
\rho_0 \circ r_R(z) & (z\in \Delta -\overline{D^{\prime\prime}} ), 
\end{cases}
$$
$\Delta$ and $\rho$ satisfy the conditions 
as in Theorem \ref{triangle}.
}

By Theorem \ref{triangle}, we obtain the following.

\thm{\label{good_triangle}
Let $(Y, q)$ be a flat surface of genus $0$
 such that $q$ has a unique zero $p_0$. 
Suppose that  $s$ is a returning geodesic on $(Y, q)$ that starts from the point $p_0$.
Then there exists a triangle $\Delta$ in $\mathbb{C}$ 
with sides $s_0, s_1, s_2$ in counterclockwise order 
and an orientation preserving  immersion 
$\rho : \Delta \to (Y, q)$ 
 satisfying the following: 
\begin{enumerate}
\item $\rho|_{\interior{\Delta}}$ is a local isometric embedding,
\item $\rho|_{\interior{\Delta}}$ contains no singular points, 
\item every vertex of $\Delta$ is mapped to $p_0$,
\item $s=\rho \circ s_0$ and 
\item\label{can_replace} $\rho \circ s_1$ is a returning geodesic.
\end{enumerate}
}

\rem{
In the condition (\ref{can_replace}) of Theorem \ref{good_triangle},
 we may replace $\rho \circ s_1$ to  $\rho \circ s_2$. 
}

\pf{
%We may assume without loss of generality that $s$ is horizontal.
Let $\vv_0$ be a unit vector that is not parallel to $s$.
By Theorem \ref{triangle}, there exists a triangle $\Delta_0$ in $\mathbb{C}$ 
with sides $s_0, s_1, s_2$  in counterclockwise order 
and an orientation preserving immersion $\rho_0 : \Delta_0 \to (Y, q)$ satisfying the following:
\begin{enumerate}[(\textrm{0-}1)]
\item $\Delta_0$ is a $(s_0, \vv_0)$-restricted triangle, 
\item $\rho_0|_{{\rm Int}(\Delta_0)}$ is a local isometric embedding,
\item $\rho_0({\rm Int}(\Delta_0))$ contains no singular points, 
\item every vertex of $\Delta_0$ is mapped to $p_0$,
\item $s=\rho_0 \circ s_0$, and 
\item $\rho_0 \circ s_i$ is a saddle connection or a returning geodesic for each $i=1, 2$.
%\item[(0-1)] $\Delta_0$ is a $(s_0, \vv)$-restricted triangle, 
%\item[(0-2)] $\rho_0|_{{\rm Int}(\Delta_0)}$ is a local isometric embedding,
%\item[(0-3)] $\rho_0({\rm Int}(\Delta_0))$ contains no singular points, 
%\item[(0-4)] every vertex of $\Delta_0$ is mapped to $p_0$,
%\item[(0-5)] $s=\rho_0 \circ s_0$, and 
%\item[(0-6)] $\rho_0 \circ s_i$ is a saddle connection or a returning geodesic for each $i=1, 2$.
\end{enumerate}
If  $s_1$  is mapped 
via $\rho_0$ to a returning geodesic, then we are done.
Now we assume that  $s_1$ is mapped 
via $\rho_0$ to a saddle connection. 
Since $(Y, q)$ is a surface of genus $0$, 
every simple closed curve on $(Y, q)$ divides $(Y, q)$.
Thus, $Y - \rho_0(\Delta_0)$ has at most two connected components.
Let $Y_1$ be the connected component of $Y- \rho_0(\Delta_0)$ 
whose boundary is $\rho_0 \circ s_1$. 
Let $q_1$ be the restriction of $q$ onto $Y_1$. 
After adjusting the parameter of $s_1$ 
so that $\left(s_1 \right)^\prime$ is constant, 
we identify $\rho_0 \circ s_1 (t)$ with $\rho_0 \circ s_1(1-t)$ 
for all $t\in [0, 1]$ on the boundary of $Y_1$.
Then $Y_1$ is a compact surface of genus $0$ 
and $q_1$ is a meromorphic quadratic differential 
such that $\rho_0 \circ s_1$ is a returning geodesic on the flat surface $(Y_1, q_1)$.
Let $\vv$ be a unit vector that is parallel to $s$.
By Theorem \ref{triangle}, there exists a triangle $\Delta_1$ in $\mathbb{C}$ 
with sides $s_1, s_3, s_4$  in counterclockwise order 
and an  orientation preserving immersion $\rho_1 : \Delta_1 \to (Y_1, q_1)$ satisfying the following:
\begin{enumerate}[(\textrm{1-}1)]
\item $\Delta_1$ is a right strongly $(s_1, \vv)$-restricted triangle, 
\item $\Delta_0 \cap \Delta_1 = s_1$, 
\item $\rho_1|_{{\rm Int}(\Delta_1)}$ is a local isometric embedding,
\item $\rho_1({\rm Int}(\Delta_1))$ contains no singular points, 
\item every vertex of $\Delta_1$ is mapped to $p_0$,
\item $\rho_0 \circ s_2=\rho_1 \circ s_2$, and 
\item $\rho_1 \circ s_i$ is a saddle connection or a returning geodesic for each $i=3, 4$
\end{enumerate}
%Moreover, we may assume the following:
%\begin{enumerate}[(\textrm{1-}1)]
%\setcounter{enumi}{6}
% and 
%%\item 
%$s_0$ and $s_4$ have a common end point
%\end{enumerate}
(see Figure \ref{process2tri}).

\begin{figure}[ht!]
\labellist
%\small
\hair 0pt
\pinlabel $s_0$  at 100 0
\pinlabel $s_1$  at 150 180
\pinlabel $s_2$  at 10 180
\pinlabel $s_3$  at 250 180
\pinlabel $s_4$  at 190 325 
\pinlabel $\Delta_0$  at 85 120
\pinlabel $\Delta_1$  at 190 260
\pinlabel $\vv$  at 350 10
 \endlabellist
\centering
 \includegraphics[scale=0.38]{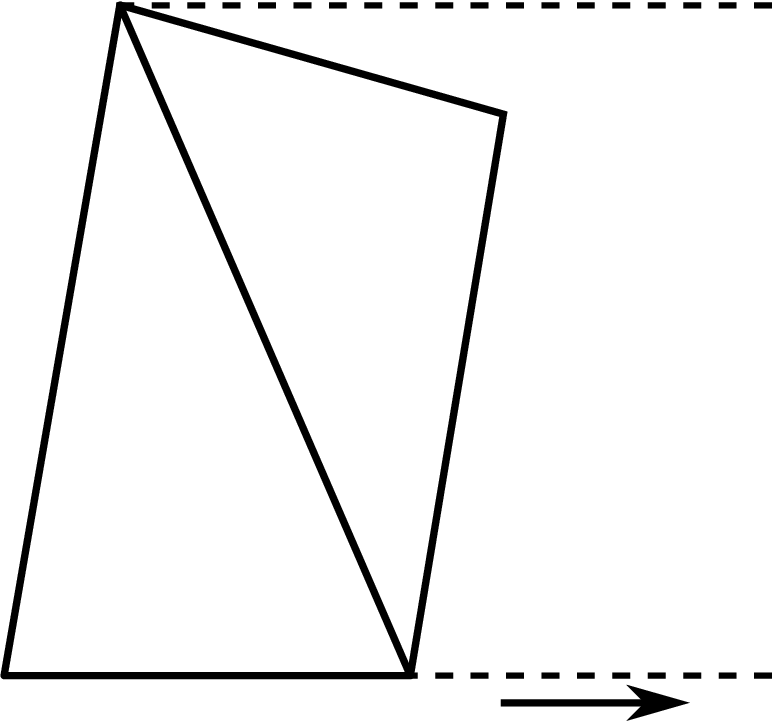} 
\caption{Triangles $\Delta_0, \Delta_1$ and dotted lines that are parallel to $\vv$.}
\label{process2tri}
\end{figure}

If $s_3$ is mapped 
via $\rho_1$ to a returning geodesic, then 
let $\Delta$ be the triangle spanned by $s_0$ and $s_3$. 
By construction, $\Delta$ is contained in $\Delta_0 \cup \Delta_1$.
Regarding $\rho_1$ as a map from $\Delta_1$ to $(Y, q)$, 
the map  $\rho : \Delta \to (Y, q)$ defined  by 
$$\rho(z)=
\begin{cases}
\rho_0(z)   & (z \in \Delta_0) \\
\rho_1 (z) & (z\in \Delta_1)
\end{cases}
$$
is well-defined. 
The triangle $\Delta$ and the map $\rho$ satisfy 
all conditions as in Theorem \ref{good_triangle}.
Hereafter we assume that  $s_3$ is not mapped 
via $\rho_1$ to a returning geodesic.
Let $Y_2$ be the connected component of $Y_1- \rho_1(\Delta_1)$ 
whose boundary is $\rho_1 \circ s_3$. 
We sew the boundary of $Y_2$ as we do for $Y_1$.
Then $Y_2$ is a compact surface of genus $0$ 
and $q_1$ induces a meromorphic quadratic differential $q_2$ on $Y_2$ 
such that $\rho_1 \circ s_3$ is a returning geodesic on the flat surface $(Y_2, q_2)$.
We apply Proposition \ref{triangle} to $\rho_1 \circ s_3$ and the vector $\vv$.
Then, there exists a triangle $\Delta_2$ in $\mathbb{C}$ 
with sides $s_3, s_5, s_6$   in counterclockwise order 
and an  orientation preserving immersion $\rho_2 : \Delta_2 \to (Y_2, q_2)$ satisfying the following:
\begin{enumerate}[(\textrm{2-}1)]
\item $\Delta_2$ is a right strongly $(s_3, \vv)$-restricted triangle, 
\item $\Delta_1 \cap \Delta_2 = s_3$, 
\item $\rho_2|_{{\rm Int}(\Delta_2)}$ is a local isometric embedding,
\item $\rho_2({\rm Int}(\Delta_2))$ contains no singular points, 
\item every vertex of $\Delta_2$ is mapped to $p_0$,
\item $\rho_1 \circ s_3=\rho_2 \circ s_3$, and 
\item $\rho_2 \circ s_i$ is a saddle connection or a returning geodesic for each $i=5, 6$
\end{enumerate}
(see Figure \ref{process3tri}). 

\begin{figure}[ht!]
\labellist
%\small
\hair 0pt
\pinlabel $s_0$  at 100 0
\pinlabel $s_1$  at 150 180
\pinlabel $s_2$  at 10 180
\pinlabel $s_3$  at 250 180
\pinlabel $s_4$  at 190 325  
\pinlabel $s_5$  at 300 140
\pinlabel $s_6$  at 290 275
\pinlabel $\Delta_0$  at 85 120
\pinlabel $\Delta_1$  at 190 260
\pinlabel $\Delta_2$  at 270 220
\pinlabel $\vv$  at 350 5
 \endlabellist
\centering
 \includegraphics[scale=0.38]{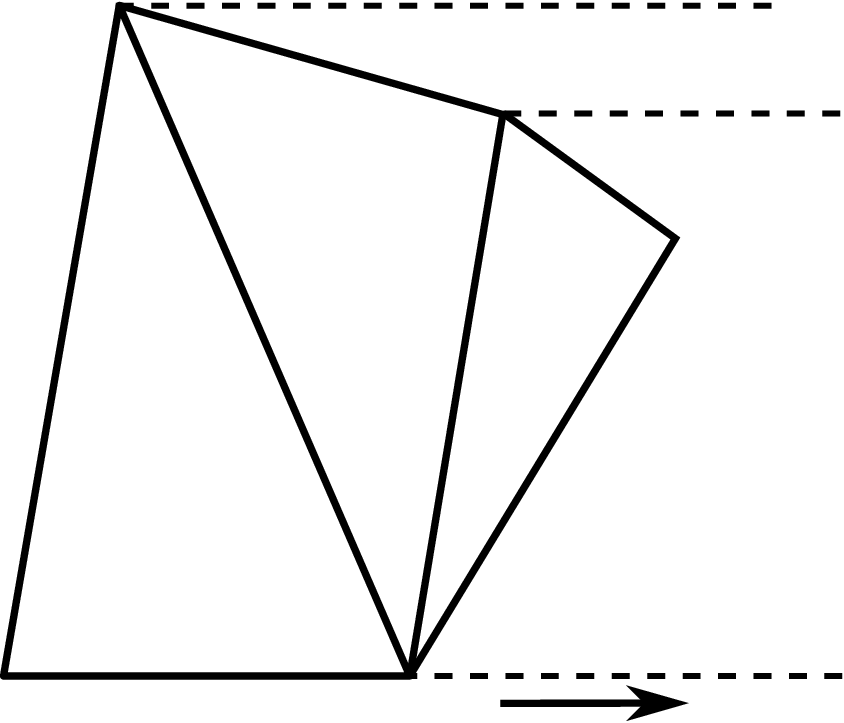} 
\caption{Triangles $\Delta_0, \Delta_1, \Delta_2$ and dotted lines that are parallel to $\vv$.}
\label{process3tri}
\end{figure}
If $s_5$ is mapped 
via $\rho_2$ to a returning geodesic, then 
let $\Delta$ be the triangle spanned by $s_0$ and $s_5$. 
By construction, $\Delta$ is contained in $\Delta_0 \cup \Delta_1 \cup \Delta_2$.
Regarding $\rho_1$ and $\rho_2$ as a map to $(Y, q)$, 
the map  $\rho : \Delta \to (Y, q)$ defined  by 
$$\rho(z)=
\begin{cases}
\rho_0(z)   & (z \in \Delta_0) \\
\rho_1 (z) & (z\in \Delta_1)\\
\rho_2 (z) & (z\in \Delta_2)
\end{cases}
$$
is well-defined. 
The triangle $\Delta$ and the map $\rho$
 satisfy all conditions as in Theorem \ref{good_triangle}.
If  $s_5$ is not mapped via $\rho_1$ to a returning geodesic, 
then we make the same argument as above 
by using the vector $\vv$. 
Since the angle around $p_0$ is finite, this process must stop by finitely many steps.
At the end of this process, we have 
%unit vectors $\vv_k$ ($k=1, \dots, n$) 
triangles $\Delta_k$ with sides $s_{2k-1}, s_{2k+1}, s_{2k+2}$ ($k=1, \dots, n$) 
and orientation preserving immersions $\rho_k : \Delta_k \to (Y_k, q_k)$  ($k=1, \dots, n$)  
satisfying the following conditions:  
\begin{enumerate}[(\textrm{3-}1)]
\item $\Delta_k$ is a right strongly $(s_{2k-1}, \vv)$-restricted triangle, 
\item $\Delta_{k-1} \cap \Delta_k = s_{2k-1}$, 
\item $\rho_k|_{{\rm Int}(\Delta_k)}$ is a local isometric embedding,
\item $\rho_k({\rm Int}(\Delta_k))$ contains no singular points, 
\item every vertex of $\Delta_k$ is mapped to $p_0$,
\item $\rho_{k-1}\circ s_{2k-1}=\rho_{k} \circ s_{2k-1}$, 
\item $\rho_k \circ s_{2k+1}$ is a saddle connection for each $i=1, \dots, n-1$, and
\item $\rho_n \circ s_{2n+1}$ is a returning geodesic
\end{enumerate}
%($k=1, \dots, n$). 
%Especially, 
% $\rho_n \circ s_{2n+2}$ is a returning geodesic 
%since  $\Delta_n$ is obtained by the final step of the process.
Let $\Delta$ be the triangle spanned by $s_0$ and $s_{2n+1}$. 
By construction, $\Delta$ is contained in $\Delta_0 \cup \Delta_1 \cup  \cdots \cup \Delta_n$.
The map  $\rho : \Delta \to (Y, q)$ defined  by 
$$\rho(z)=
\begin{cases}
\rho_0(z)   & (z \in \Delta_0) \\
\rho_1 (z) & (z\in \Delta_1)\\
 & \vdots \\
\rho_n (z) & (z\in \Delta_n)
\end{cases}
$$
is well-defined. 
Now $\Delta$ and $\rho$ satisfy all conditions as in Theorem \ref{good_triangle}.
}

\subsection{Proof of the inequality (\ref{lower_bound})}
In this subsection, we prove Theorem \ref{main}.
The upper estimates of $N(X, \omega)$ are given by Proposition \ref{upper_bound_proof}. 
We give the lower estimates in Proposition \ref{lower_bound_proof}.

\lem{\label{choose_delta1}
Let $g\geq 2$. 
Let $(Y, q)$ be a flat surface of genus $0$ 
that has a unique zero $p_0$ of order $2g-3$. 
Fix a simple pole $p_1$ of $q$ and 
a returning geodesic $s_0$ on $(Y, q)$ 
that passes through $p_1$.
Then there exist saddle connections $\delta_1, \dots, \delta_g$ 
satisfying the following conditions:
\begin{itemize}
\item the endpoints of  $\delta_1, \dots, \delta_g$
 are poles of $q$ and 
\item  $\delta_1, \dots, \delta_g$, and $s_0$
are disjoint to each other.
\end{itemize}
}

%To prove Lemma\ref{choose_delta1}, 
%we give the following definition.
%
%
%\defi{
%Let $(Y, q)$ be a flat surface of genus $0$ 
%that has a unique zero 
%and an odd number of poles. 
%Fix a returning geodesic $s_0$
%A saddle connection $s$ is odd  
%if 
%}

\pf{We prove it by induction on $g$.
If $g=2$ then $p_0$ is a unique zero of $q$ 
of order $1$. 
By Theorem \ref{good_triangle}, 
there exists 
a triangle $\Delta_0$ in $\mathbb{C}$ 
with sides $s_0, s_1, s_2$ in counterclockwise order 
and an orientation preserving  immersion 
$\rho_0 : \Delta_0 \to (Y, q)$ 
 satisfying the 
following conditions: 
\begin{enumerate}[(\textrm{0-}1)]
\item $\rho_0|_{\interior{\Delta_0}}$ is a local isometric embedding,
\item $\rho_0|_{\interior{\Delta_0}}$ contains no singular points, 
\item every vertex of $\Delta_0$ is mapped to $p_0$,
\item $s=\rho \circ s_0$, and 
\item $\rho _0\circ s_1$ is a returning geodesic.
\end{enumerate}
% same conditons 
% as in Theorem \ref{good_triangle}. 
 Let $P_0$ be the vertex of $\Delta_0$ that is opposite to $s_1$ 
 and $P_1$ the midpoint of $s_1$. 
The surface $Y_1=Y- \rho_0(\Delta_0)$ 
is of genus $0$ and  
has a boundary  $\rho_0 \circ s_2$. 
Let $q_1$ be the restriction of $q$ onto $Y_1$. 
After adjusting the parameter of $s_2$ 
so that $\left(s_2 \right)^\prime$ is constant, 
we identify $\rho_0 \circ s_2 (t)$ with $\rho_0 \circ s_2(1-t)$ 
for all $t\in [0, 1]$ on the boundary of $Y_1$.
Then $Y_1$ is a compact surface of genus $0$ 
and $q_1$ is a meromorphic quadratic differential 
such that $\rho_0 \circ s_2$ is 
a returning geodesic on the flat surface $(Y_1, q_1)$.
Let $\vv$ be a unit vector  that is parallel to the segment $P_0P_1$.
By Theorem \ref{triangle},  
 there exists  
a triangle $\Delta_1$ in $\mathbb{C}$ 
with sides $s_2, s_3, s_4$
and an  orientation preserving  immersion 
$\rho_1 : \Delta_1 \to (Y_1, q_1)$ 
 satisfying the following conditions (see Figure \ref{find_delta_g=2-1}): 
\begin{enumerate}[(\textrm{1-}1)]
\item $\Delta_1$ is a left  strongly $(s_2, \vv)$-restricted triangle, 
\item $\Delta_0 \cap \Delta_1= s_2$,
%\item $s_2$ is parallel to $s$ and 
\item $\rho_1|_{{\rm Int}(\Delta_1)}$ is a local isometric embedding,
\item $\rho_1({\rm Int}(\Delta_1))$ contains no singular points, 
\item every vertex of $\Delta_1$ is mapped to $p_0$,
\item $\rho_1 \circ s_2=\rho_0 \circ s_2$, and 
\item $\rho_1 \circ s_i$ is a saddle connection or a returning geodesic for each $i=3,4$.
%\item If $s_i$ is parallel to $\vv$ for some $i=1,2$, then $\rho \circ s_i$ is a returning geodesic.
\end{enumerate}

\begin{figure}[h]
\labellist
\hair 0pt
\pinlabel $\Delta_0$  at    270 90
\pinlabel $\Delta_1$  at    150 120
\pinlabel $s_0$  at    250 20
\pinlabel $s_1$  at    350 120
\pinlabel $\vv$  at   140 200
\pinlabel $s_2$  at    270 170
\pinlabel $s_3$  at    190 160
\pinlabel $s_4$  at     30 55
\pinlabel $P_0$  at    110 20
\pinlabel $P_1$  at    370 170
\endlabellist
\centering
\includegraphics[scale=0.38]{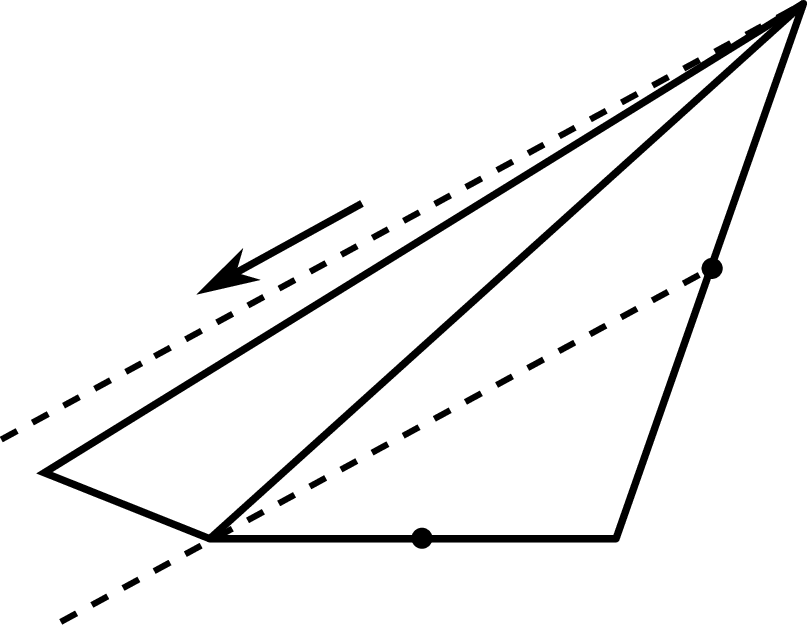} 
\caption{The triangles $\Delta_0$ and $\Delta_1$}
\label{find_delta_g=2-1}
\end{figure}
If  both 
$\rho_1 \circ s_3$ and $\rho_1 \circ s_4$ 
are returning geodesics,  
we have
$Y=\rho_0(\Delta_0) \cup \rho_1(\Delta_1)$.
This contradicts the assumption that 
$q$ has a zero $p_0$.
Thus one of 
$\rho_1 \circ s_3$ and $\rho_1 \circ s_4$  
is a saddle connection. 
Here, we assume that $\rho_1 \circ s_3$ 
is a saddle connection.
For the case where $\rho_1 \circ s_4$ 
is a saddle connection, 
we can apply the same argument. 
The surface $Y_2=Y_1- \rho_1(\Delta_1)$ 
is of genus $0$ and  
has a boundary  $\rho_1 \circ s_3$. 
Let $q_2$ be the restriction of $q_1$ onto $Y_2$. 
We sew the boundary of $Y_2$ as we do for $Y_1$.
Then $Y_2$ is a compact surface of genus $0$ 
and $q_2$ is a meromorphic quadratic differential 
such that $\rho_1 \circ s_3$ is 
a returning geodesic on the flat surface $(Y_2, q_2)$.
Let $\vv^\prime$ be a unit vector that is not parallel to $s_3$. 
By Theorem \ref{triangle},  
 there exists  
a triangle $\Delta_2$ in $\mathbb{C}$ 
with sides $s_3, s_5, s_6$
and an  orientation preserving  immersion 
$\rho_2 : \Delta_2 \to (Y_2, q_2)$ 
 satisfying the following conditions: 
\begin{enumerate}[(\textrm{2-}1)]
\item $\Delta_2$ is a left  strongly $(s_3, \vv)$-restricted triangle, 
\item $\Delta_1 \cap \Delta_2= s_3$,
%\item $s_2$ is parallel to $s$ and 
\item $\rho_2|_{{\rm Int}(\Delta_2)}$ is a local isometric embedding,
\item $\rho_2({\rm Int}(\Delta_2))$ contains no singular points, 
\item every vertex of $\Delta_2$ is mapped to $p_0$,
\item $\rho_2 \circ s_3=\rho_1 \circ s_3$, and 
\item $\rho_2 \circ s_i$ is a saddle connection or a returning geodesic for each $i=5,6$.
%\item If $s_i$ is parallel to $\vv$ for some $i=1,2$, then $\rho \circ s_i$ is a returning geodesic.
\end{enumerate}
Since the angle around $p_0$ in $(Y,  q)$ is $3\pi$
and all vertices of $\Delta_j$  are mapped 
to $p_0$ by $\rho_i$ for each $i=0,1, 2$,  
we have 
$Y=\rho_0(\Delta_0) \cup \rho_1(\Delta_1)\cup \rho_2(\Delta_2)$. 
Thus 
$\rho_1 \circ s_4$, 
$\rho_2 \circ s_5$, and 
$\rho_2 \circ s_6$  
are returning geodesics. 
Set $\Delta=\Delta_0\cup \Delta_1 \cup \Delta_2$ 
define the map  $\rho : \Delta \to (Y, q)$ by 
$$\rho(z)=
\begin{cases}
\rho_0(z)   & (z \in \Delta_0) \\
\rho_1 (z) & (z\in \Delta_1)\\
\rho_2 (z) & (z\in \Delta_2).
\end{cases}
$$
Let $P_i$ be the midpoint of the side $s_i$ 
for $i=4, 5,6$ 
 (see Figure \ref{find_delta_g=2-2}). 
Then $\delta_1=\rho(P_1P_4)$ 
and $\delta_2=\rho(P_5P_6)$  
are saddle connections satisfying all conditions 
as in Lemma \ref{choose_delta1}.

\begin{figure}[h]
\labellist
\hair 0pt
\pinlabel $\Delta_0$  at    270 90
\pinlabel $\Delta_1$  at    150 120
\pinlabel $s_0$  at    250 20
\pinlabel $s_1$  at    350 120
\pinlabel $\Delta_2$  at   140 200
\pinlabel $s_2$  at    270 170
\pinlabel $s_3$  at    190 160
\pinlabel $s_4$  at    30 55
\pinlabel $P_0$  at    110 20
\pinlabel $P_1$  at    370 170
\pinlabel $P_4$  at    60 30
\pinlabel $P_5$  at    200 300
\pinlabel $P_6$  at    0 160
\pinlabel $s_5$  at    140 280
\pinlabel $s_6$  at    15 200
\endlabellist
\centering
\includegraphics[scale=0.38]{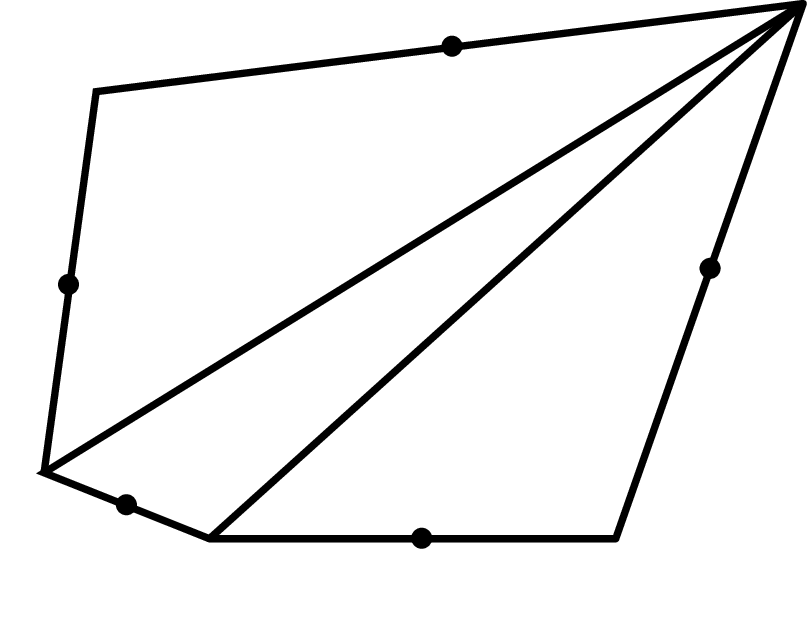} 
\caption{The triangles $\Delta_0$, $\Delta_1$, and $\Delta_2$}
\label{find_delta_g=2-2}
\end{figure}

Next, we assume that Lemma \ref{choose_delta1} 
is true for all flat surfaces of genus $0$ 
that have a unique zero $p_0$ of odd order  less than $2g-3$. 
Let $(Y, q)$ be a flat surface of genus $0$ 
that has a unique zero $p_0$ of order $2g-3$. 
Let $s_0$ be a returning geodesic on $(Y, q)$. 
We can construct triangles $\Delta_1$, $\Delta_2$, 
immersions $\rho_1: \Delta_1\to (Y, q)$ 
and $\rho_2: \Delta_2 \to (Y, q)$   
that are the same as in the case of $g=2$. 
We use the same labels for the sides and vertices.
If $s_i \in \left\{ s_3, s_4 \right\} $ corresponds to a returning geodesic, 
then the segment connecting the midpoint of $s_i$ and $P_1$ 
corresponds to a saddle connection connecting 
poles of $(Y, q)$. 
We denote the saddle connection by $\delta_1$.
Since $(Y, q)$ has $2g-3$ poles, 
the side $s_{7-i}$ is a saddle connection. 
In the same way as in the case of $g=2$, 
We can regard $Y_2=Y-\left(\rho_1(\Delta_1)\cup \rho_2(\Delta_2)\right)$ 
as a compact surface of genus $0$ without boundary and 
$q_2=q|_{Y_2}$ 
as a meromorphic quadratic differential  
with $2g-5$ poles.  
Moreover, $\rho(s_{7-i})$ is a returning geodesic on $(Y_2, q_2)$.
By assumption, we can find 
saddle connections $\delta_2, \dots, \delta_g$ 
in $(Y_2, q_2)$ 
satisfying all conditions as in Lemma \ref{choose_delta1}.
Regarding them as saddle connections in $(Y, q)$,  
$\delta_1, \dots, \delta_g$ are saddle connections in $(Y, q)$ 
satisfying all conditions as in Lemma \ref{choose_delta1}.
Next, assume that both $s_3$ and $s_4$ are saddle connections. 
Then $Y-\left(\rho_1(\Delta_1)\cup \rho_2(\Delta_2)\right)$  
has two connected components.
% $Y^\prime$ and $Y^{\prime\prime}$.
Since $(Y, q)$ has an odd number of poles, 
one of the components 
%we can assume that $Y^\prime$
contains an odd number of poles (before sewing its boundary) 
and the other 
%$Y^{\prime\prime}$  
contains an even number of poles.
Hereafter, we say a saddle connection $s^\prime$ 
is odd if 
the connected component of $Y-s^\prime$ 
that does not contain $s_0$ 
has an odd number of poles.
We delete the labels $s_3$ and $s_4$ 
and construct triangles $\Delta_2, \Delta_3, \dots $ 
and immersions $\rho_i : \Delta_i \to (Y_i, \rho_i)$ ($i=2, 3, \dots$) 
by the following process  (see Figure \ref{process}):
\begin{enumerate}
\item Set $i=1$. 
\item\label{go_back} 
Let
$s_{i+2}$ be the side of  $\Delta_{i}$  
such that  $\rho_{i}(s_{i+2})$  is  odd and 
$P_{i+1}$ the vertex of  $\Delta_{i}$ 
that is opposite to $s_{i+1}$. 
Let $\vv_{i}$ 
be a unit vector that is parallel to the segment $P_1P_{i+1}$.   
Set  
$\displaystyle Y_{i+1}=Y-\bigcup_{k=0}^i \rho_i(\Delta_i) $ 
and 
regard it as a closed surface of genus $0$ 
without boundary by sewing $\rho_i(s_{i+2})$. 
Set $q_{i+1}=q|Y_{i+1}$.
\item 
By applying Theorem \ref{triangle}, let 
 $\Delta_{i+1}$ be 
 a triangle in  $\mathbb{C}$
with $s_{i+2}$ as a side 
and  
$\rho_{i+1} : \Delta_{i+1 }\to (Y_{i+1}, q_{i+1})$ 
orientation preserving  immersion 
satisfying the following conditions : 
\begin{enumerate}[(i)]
\item $\Delta_{i+1}$ is a left  strongly $(s_{i+2}, \vv_i)$-restricted triangle, 
\item $\Delta_{i+1} \cap \Delta_i= s_{i+2}$,
%\item $s_2$ is parallel to $s$ and 
\item $\rho_{i+1}|_{{\rm Int}(\Delta_{i+1})}$ is a local isometric embedding,
\item $\rho_{i+1}({\rm Int}(\Delta_{i+1}))$ contains no singular points, 
\item every vertex of $\Delta_{i+1}$ is mapped to $p_0$,
\item $\rho_{i+1} \circ s_{i+2}=\rho_i \circ s_{i+2}$, and 
\item other sides of $\Delta_{i+1}$ than $s_{i+2}$ are mapped via $\rho_{i+1}$ 
saddle connections or returning geodesics.
%\item If $s_i$ is parallel to $\vv$ for some $i=1,2$, then $\rho \circ s_i$ is a returning geodesic.
\end{enumerate}
\item If other sides of $\Delta_{i+1}$ than $s_{i+2}$ are both 
mapped via $\rho_{i+1}$ 
saddle connections, 
then we replace $i$ to $i+1$ and repeat the process from (\ref{go_back}).
If not, we stop this process.
\end{enumerate} 
\begin{figure}[h]
\labellist
\hair 0pt
\pinlabel $\Delta_0$  at    470 90
\pinlabel $\Delta_1$  at    370 110
\pinlabel $\Delta_2$  at    290 120
\pinlabel $\Delta_3$  at    220 160
\pinlabel $\Delta_{i+1}$  at  50 160
\pinlabel $s_0$  at   450 -10
\pinlabel $s_1$  at    540 100
\pinlabel $s_2$  at    420 100
\pinlabel $s_3$  at    280 56
\pinlabel $s_4$  at    240 105
\pinlabel $s_5$  at   180 130
\pinlabel $s_{i+2}$  at  90 110
\pinlabel $P_0$  at    300 -20
\pinlabel $P_1$  at    560 130
\pinlabel $P_2$  at    220 10
\pinlabel $P_3$  at    135 50
\pinlabel $P_4$  at    200 210
\endlabellist
\centering
\includegraphics[scale=0.38]{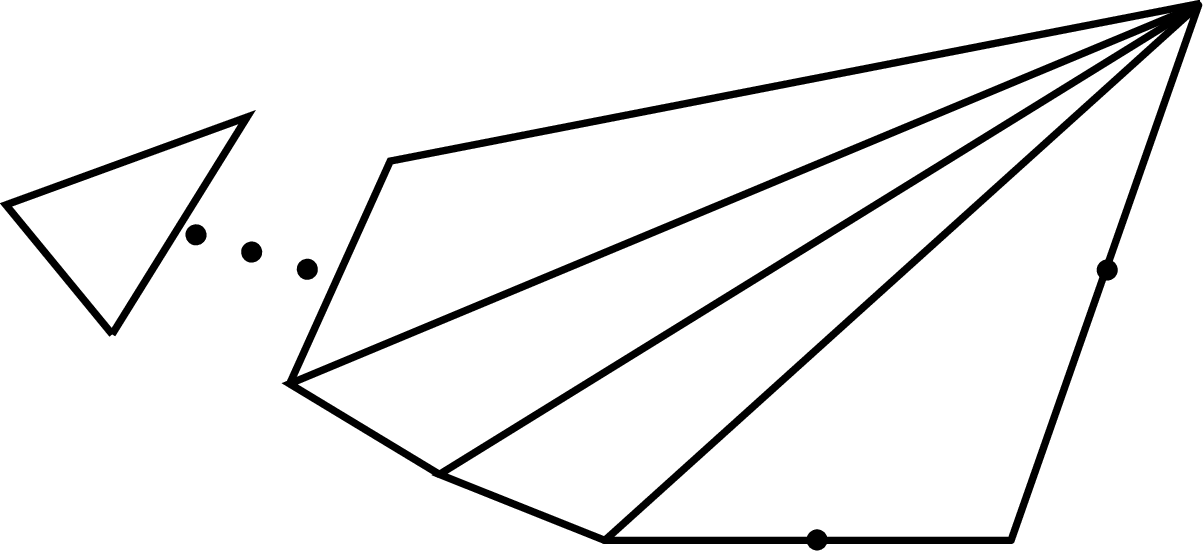} 
\caption{The process to construct triangles  $\Delta_2, \Delta_3, \dots $}
\label{process}
\end{figure}
Note that 
the segment $P_1P_{i+1} $ ($i=2,3, \dots$) 
is contained in 
$\displaystyle \bigcup_{k=0}^i \Delta_i $ 
. 
This process stops by finitely many steps 
since the angle around $p_0$ is finite.
Assume that we obtain 
triangles $\Delta_2, \Delta_3, \dots, \Delta_n$ 
and immersions $\rho_i : \Delta_i \to (Y_i, \rho_i)$ 
($i=2, 3, \dots, n$) satisfying the above conditions 
as the result of the process. 
Then $\Delta_n$ has a side $s_{n+1}$  
that is mapped via $\rho_n$ a returning geodesic. 
Let $P_{n+1}$ be the midpoint of $s_{n+1}$. 
By construction, the segment $P_1 P_{n+1}$ 
is in 
$\displaystyle \bigcup_{k=0}^n \Delta_k $. 
The union 
$\delta_1= \displaystyle \bigcup_{k=0}^n \rho_k\left(P_1 P_{n+1}\cap \Delta_k\right) $ 
is a saddle connection in $(Y, q)$ 
that is disjoint from $s_0$.
Let $W_1, W_2, \dots, W_n$ be  the 
connected  components of  
$\displaystyle Y- \bigcup_{k=0}^n \rho_k(\Delta_k)$. 
We set $q_i=q|_{W_i}$ for each $i=1, 2, \dots, n$. 
The flat surface $(W_i, q_i)$ is of genus $0$ and has  a boundary.
Denote by $2g_i$ the number of poles of $(W_i, q_i)$.
Sewing the boundary,  
$(W_i, q_i)$ 
is a flat surface that has $2g_i+1$ poles  
and the boundary is now a returning geodesic, say $s^\prime_i$.
By the assumption of the induction, 
we can find saddle connections $\delta^i_{1}, \delta^i_{2}, \dots, \delta^i_{g_i}$ 
satisfying the following:
\begin{itemize}
\item the endpoints of 
 $\delta^i_{1}, \delta^i_{2}, \dots, \delta^i_{g_i}$ 
are poles of $q_i$ and 
\item  $\delta^i_{1}, \delta^i_{2}, \dots, \delta^i_{g_i}, s^\prime_i$ 
are disjoint to each other.
\end{itemize}
Since $\displaystyle \sum_{i=1}^n g_i=g-1 $, 
the collection of curves 
\begin{align*}
\left\{\delta_1 \right\}
\cup 
\left\{\delta^i_{j} : 
i\in\left\{1,2, \dots, n\right\} , 
j\in\left\{1, 2, \dots, g_i \right\} \right \}
\end{align*}
contains $g$ saddle connections. 
The saddle connections satisfy 
all  conditions in Lemma\ref{choose_delta1}.
}

Next we prove the case 
that corresponds to the case of $\calhhyp(g-1, g-1)$.

\lem{\label{choose_delta2} 
Let $g\geq 2$. 
Let $(Y, q)$ be a flat surface of genus $0$ 
that has a unique zero $p_0$ of order $2g-2$. 
Then there exist saddle connections $\delta_1, \dots, \delta_{g+1}$ 
such that each of whose end points are poles  
and they are disjoint to each other.
}

\pf{
Given a returning geodesic $s$ on $(Y, q)$. 
By Theorem \ref{good_triangle}, 
there exist a triangle $\Delta$ in $\mathbb{C}$ 
with sides $s_0, s_1, s_2$ 
in counterclockwise order 
and an orientation preserving immersion 
$\rho: \Delta \to (Y, q)$ 
that satisfy the same conditions as in Theorem \ref{good_triangle}.
We denote by $\delta_{g+1}$ be the saddle connection in $\rho(\Delta)$ 
that connects two poles.  
We set $Y^\prime=Y-\rho(\Delta)$ and 
$q^\prime=q|_{Y^\prime}$. 
Then the surface $(Y^\prime, q^\prime)$ 
is of genus $0$ and  has a boundary. 
Sewing the boundary, 
$(Y^\prime, q^\prime)$  
is a flat surface that has $2g+1$ poles. 
By Lemma \ref{choose_delta1}, 
there exist 
 saddle connections $\delta_1, \dots, \delta_g$ 
satisfying the following conditions:
\begin{itemize}
\item the endpoints of  $\delta_1, \dots, \delta_g$
 are poles of $q$ and 
\item  $\delta_1, \dots, \delta_g$, and $s_0$
are disjoint to each other.
\end{itemize}
Regarding them as saddle connections on $(Y, q)$, the 
saddle connections $\delta_1, \dots, \delta_{g+1}$ 
satisfy all conditions as desired.
}

By Lemma \ref{choose_delta1} 
and Lemma \ref{choose_delta2}, 
we have the following.
  
\prop{\label{lower_bound_proof}
Let $g\geq 2$. 
%We have the following.
\begin{enumerate}
\item If $ (X, \omega) \in \calhhyp (2g-2)$, then 
there exist disjoint regular closed geodesics $\gamma_1, \dots, \gamma_g$ on $(X, \omega)$ 
that are invariant under $\tau$ and are not homotopic to each other. 
Moreover, we have $N(X, \omega) \geq g$.
\item \label{even_case} If $(X, \omega) \in \calhhyp(g-1, g-1)$, then 
there exist disjoint regular closed geodesics $\gamma_1, \dots, \gamma_{g+1}$ on $(X, \omega) $
that are invariant under $\tau$ and are not homotopic to each other. 
Moreover, we have $N(X, \omega) \geq g+1$.
\end{enumerate} 
}

\pf{
Assume that $ (X, \omega) \in \calhhyp (2g-2)$ 
and $\tau$ is the hyperelliptic involution of $(X, \omega)$.
Let $(Y, q)$ be a flat surface of genus $0$ 
that is the quotient  $(X, \omega)$ by $\tau$.
Then $(Y, q)$ has a unique zero $p_0$ of order $2g-3$. 
By Lemma \ref{choose_delta1}, 
there exist disjoint saddle connections 
$\delta_1, \dots, \delta_g$ 
each of which connects distinct poles. 
Let $\gamma_i$ be the regular closed geodesic in $(X, \omega)$ 
that is projected to $\delta_i$ 
by the natural projection from $(X, \omega)$ to $(Y, q)$ for each $i=1, \dots, g.$ 
Then the regular closed geodesics 
$\gamma_1 \dots, \gamma_g$ 
are disjoint to each other. 
Therefore, we have $N(X, \omega)\geq g$. 
The proof of the statement of (\ref{even_case}) 
is also done in the same way.
%and $2g+1$ simple poles.
%Let $s$ be a returning geodesic on $(Y, q)$. 
%By Corollary \ref{quadrilateral}, 
%there exist a quadrilateral $Q_1$ 
%with sides $s_0, s_1, s_2, s_3$ in counterclockwise order 
%and $\rho_1 : Q_1 \to \mathbb{C}$ 
%satisfying the same conditions as in 
}

\section{Classification theorem for translation surfaces in $\calhhyp(4)$}

The classification theorem for compact surfaces claims 
that every  oriented closed surface  of genus $g\geq 1$
is topologically constructed from a sphere  
by gluing $g$ cylinders.
In this section, we prove the same kind of classification theorem for
almost all translation surfaces in $\calhhyp(4)$ 
with respect to their  structures as translation surfaces.

\defi{[Maximal cylinder and simple cylinder]
Let $(X, \omega)$ be a translation surface and 
$\gamma$ a regular closed geodesic on $(X, \omega)$.
A \emph{maximal cylinder} for  $\gamma$ is the union of all regular closed geodesics 
%on a translation surface 
that are homotopic to $\gamma$ (see Proposition \ref{geod_rep}). 
A \emph{simple cylinder} is a maximal cylinder
each of whose boundary components is only one saddle connection. 
}

\thm{\label{classification}
Let $(X, \omega) \in \calhhyp(4)$ 
and $\tau$ the hyperelliptic  involution of $(X, \omega)$.
If $(X, \omega) \not \in \gltr \cdot \str$, 
then there exist disjoint simple cylinders $C_1, C_2, C_3$ each of which  is 
invariant under $\tau $. 
If $(X, \omega) \in \gltr \cdot \str$, 
then $(X, \omega)$ has at most two  disjoint simple cylinders. 
}
 
 As corollaries of this theorem, we have the following.

\begin{figure}[h]
\labellist
\hair 0pt
\pinlabel $p_0$  at   66 242
\pinlabel $p_1$  at   45  215 
\pinlabel $p_2$  at  56 153 
\pinlabel $p_3$  at   145 138
\pinlabel $p_4$  at  210 182
\pinlabel $p_5$  at  185 240 
\pinlabel $p_6$  at  105 258
\pinlabel $C_1$  at  425 250
\pinlabel $C_2$  at  400 135 
\pinlabel $C_3$  at  380 45
\pinlabel $s_1$  at  68 184
\pinlabel $s_2$  at  100 160 
\pinlabel $s_3$  at  142 167
\endlabellist
\centering
\includegraphics[scale=0.55]{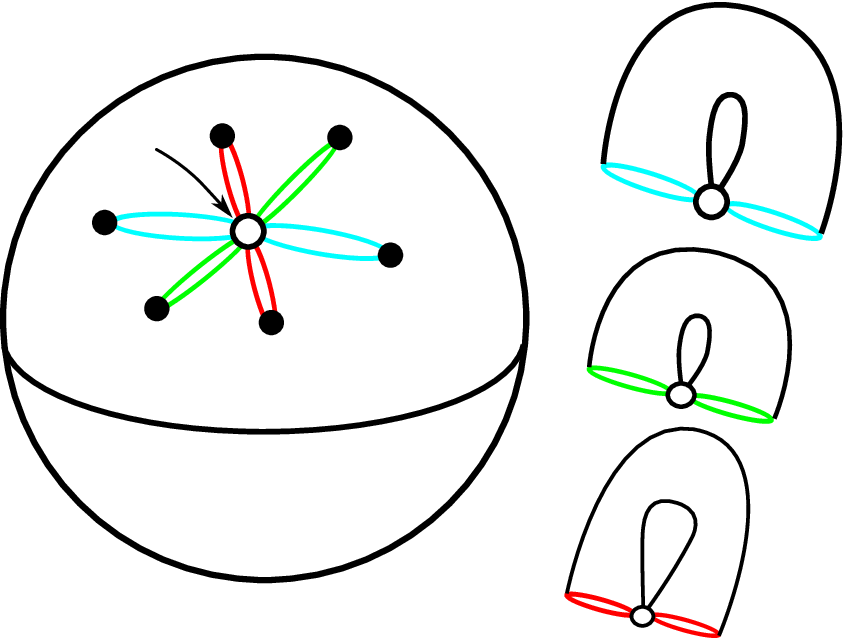} 
\caption{The flat surface $(X_0, q)$ and cylinders $C_1, C_2$, and $C_3$.}  
\label{classification3}
\end{figure}

 \cor{[Classification theorem for translation surfaces in $\calhhyp(4)$]\label{classification_cor}
 Let $(X, \omega) \in \calhhyp(4) \setminus \gltr \cdot \str$. 
Then, there exists a flat surface $(X_0, q)$ of genus $0$ 
such that  
 $q$ 
has a zero $p_0$ of order $2$ and six simple poles $p_1, p_2, \dots, p_6$, 
an involution $\tau_0 : X_0 \to X_0$, 
saddle connections $s_1, s_2, s_3$ and 
Euclidean cylinders $C_1, C_2, C_3$
satisfying the following (see Figure \ref{classification3}):
\begin{enumerate}
\item $\tau_0^\ast q=q$, 
\item $s_i$ is a saddle connection of $(X_0, q)$ connecting $p_0$ and $p_i$ for $i=1, 2, 3$, 
\item $\tau_0(s_i)$ connects $p_0$ and $p_{i+3}$ for $i=1, 2, 3$, 
\item two of each $s_1, s_2, s_3, \tau_0(s_1), \tau_0(s_2)$, and $\tau_0(s_3)$ 
         intersect  only at $p_0$,
\item the circumference of $C_i$ equals $2|s_i|$ for $i=1,2,3$, and
\item  $(X, \omega)$ is obtained by cutting $(X_0, q)$ along the saddle connections 
$s_1$, $s_2$, $s_3$, $\tau_0(s_1)$, $\tau_0(s_2)$, $\tau_0(s_3)$ 
and gluing $C_i$ to the slits $s_i$ and $\tau_0(s_i)$ for all $i=1, 2, 3$.
\end{enumerate}
 }

\pf{
Let $(X, \omega) \in \calhhyp(4) \setminus \gltr \cdot \str$.  
By Theorem \ref{classification}, $(X, \omega) $ has 
disjoint simple cylinders $C_1, C_2, C_3$ each of which  is 
invariant under $\tau $. 
Removing the cylinders from $(X, \omega)$ and sewing the $6$ boundaries respectively,  
the resulting surface $X_0$ is a Riemann surface of genus $0$ 
and has a quadratic differential $q$ induced by $\omega^2$. 
The hyperelliptic involution $\tau$ of $(X, \omega)$ induces 
an involution $\tau_0: X_0 \to X_0$. 
Let $s_i$ be the sewed boundary of $(X_0, q)$. 
The boundary component $s_i$ is also  a boundary component of $C_i$ 
for each $i=1, 2, 3$. 
For each $i=1, 2, 3$, the other boundary component of $C_i$ is $\tau(s_i)$ 
since $C_i$ is invariant under $\tau$ and $\tau^\ast \omega=-\omega$.
Therefore, $\tau_0$ and $s_1, s_2, s_3$ satisfy the conditions as in the claim.
}

%Let $(Y, q)$ be a flat surface constructed from $(X, \omega)$
%by removing  the cylinders $C_1, C_2, C_3$  
%and sewing the $6$ boundaries respectively. 
%Let $p_0$ be the unique zero of $\omega$ 
%and $p_1$ the point of $(Y, q)$ corresponding to $p_0$. 
%Since $C_1, C_2, C_3$ are simple cylinders 
%and the angle ground $p_0$ on $(X, \omega)$ is $10\pi$, 
%the angle around $p_1$ is $4\pi$. 
%This implies that $p_1$ is a zero of $q$ of order $2$. 

 \cor{\label{hex_para_cor}
 Let $(X, \omega) \in \calhhyp(4) \setminus \gltr \cdot \str$. 
 Then, $(X, \omega)$ is constructed from a center-symmetric hexagon $H$ 
 and $3$ parallelograms $P_1, P_2, P_3$ by gluing them as in Figure \ref{hex_para}.
 \begin{figure}[h]
\labellist
\hair 0pt
\pinlabel $H$  at  189 88
\pinlabel $P_1$ at 320 152 
\pinlabel $P_2$ at 177 203
\pinlabel $P_3$ at 56 112
\pinlabel $a_1$ at  358 108
\pinlabel $a_1$ at  293 194
\pinlabel $a_2$ at 252 209
\pinlabel $a_2$ at 90 209
\pinlabel $a_3$ at 70 180
\pinlabel $a_3$ at 42 40 
\pinlabel $b_1$ at 98 25
\pinlabel $b_1$ at 378 169
\pinlabel $b_2$ at 164 248  
\pinlabel $b_2$ at 202 -10
\pinlabel $b_3$ at 8 116
\pinlabel $b_3$ at 300 52
\endlabellist
\centering
\includegraphics[scale=0.55]{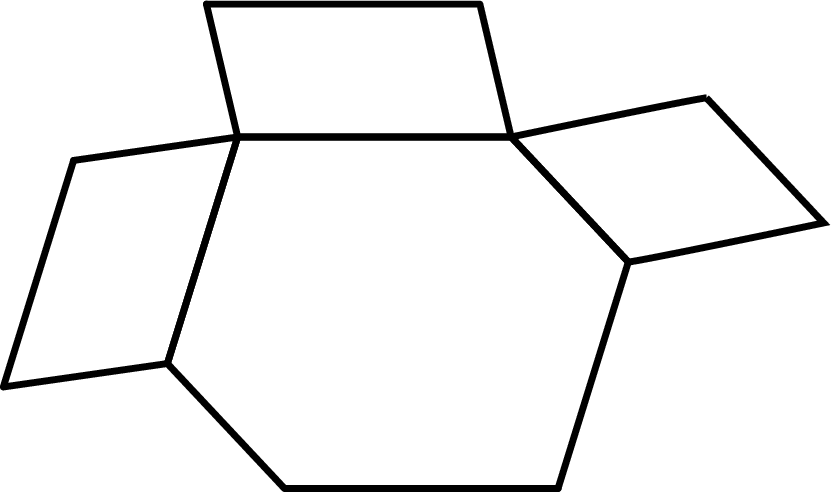} 
\caption{Construction of $(X, \omega) $ from a center-symmetric  hexagon 
and $3$ parallelograms.}
\label{hex_para}
\end{figure}
 }
 
\pf{
Let  $(X, \omega) \in \calhhyp(4) \setminus \gltr \cdot \str$. 
By Corollary \ref{classification_cor}, 
the hyperelliptic translation surface $(X, \omega)$ is constructed from 
a flat surface  $(X_0, q)$ that has a unique zero $p_0$ and 
an involution $\tau_0$ 
by gluing Euclidean cylinders $C_1, C_2, C_3$ 
along slits $s_1, s_2, s_3, \tau_0(s_1), \tau_0(s_2)$, and $\tau_0(s_3)$. 
We cut $X_0$ along the slits.
Since $\tau_0$ is an involution, the resulting surface is a center-symmetric hexagon.
We also cut each cylinder $C_i$ ($i=1, 2, 3$) 
along a segment connecting the points corresponding  to $p_0$ on its boundary components.
Then, the resulting surfaces are parallelograms.
}
 
Hereafter, 
we show  Theorem \ref{classification}. 
First, we prove the last part of Theorem \ref{classification}.

\lem{\label{only_two}
If $(X, \omega) \in \gltr \cdot \str$, 
then $(X, \omega)$ has at most two  disjoint simple cylinders. 
}

\pf{
We may assume that $(X, \omega)=\str$. 
Then, the Veech group $\Gamma(\str)$ is a lattice in $\sltr$ and 
the quotient $\mathbb{H}/\Gamma(\str)$ has two punctures (see Example \ref{exam_st5}). 
Let $p_1$ and $p_2$ be fixed points of the hyperelliptic involution that are 
the images of the midpoints of sides with labels $a$ and $3$ as in Figure \ref{str}, respectively. 
%By Proposition \ref{peripara},  we have 
%${\rm Peri}(X, \omega)/\Gamma(X, \omega)=\left\{ \left[0 \right], \left[\frac{\pi}{4} \right] \right\}$.
The cylinder decomposition for direction $\frac{\pi}{4}$ has only one cylinder that is not a simple cylinder. 
The cylinder decomposition for direction $0$ has three cylinders. 
Only one of the cylinders is a simple cylinder. 
Let $C_0$ denote this simple cylinder.  
Let $C$ be a simple cylinder on $(X, \omega)$.
By Theorem \ref{veech_dic},  
the cylinder $C$ is contained in a cylinder decomposition 
for some direction $\theta$. 
Since ${\rm Peri}(X, \omega)/\Gamma(X, \omega)=\left\{ \left[0 \right], \left[\frac{\pi}{4} \right] \right\}$,  
there exists $A\in \Gamma(X, \omega)$ 
that maps the direction $0$ to $\theta$.
By Proposition \ref{affine_homeo} and Proposition \ref{preserve_fp}, the cylinder $C$ must contain $p_1$ or $p_2$.
This means that $(X, \omega)$ has at most two  disjoint simple cylinders. 
}

Next, we prove the first claim of Theorem \ref{classification}.
Suppose that 
$(X, \omega) \in \calhhyp(4)$ 
and $\tau$ is the hyperelliptic involution of $(X, \omega)$.
We assume that $(X, \omega)$ has no three
disjoint simple cylinders each of which  is 
invariant under $\tau $. 
We will show that $(X, \omega) \in \gltr \cdot \str$. 
To do this, we set $Y=X/\left<\tau \right>$. 
Let $q$ be the holomorphic quadratic differential on $Y$ 
induced by $\omega^2$ via the natural projection 
$\varphi: X \to Y$. 
Then  the genus of $Y$ is $0$ 
and the quadratic differential $q$ has a unique zero $p_0$ of order $3$ 
and $7$ simple poles. 
As we do in Section \ref{proof_of_main_thm}, 
we construct the union of some triangles 
whose sides correspond to saddle connections 
or returning geodesics  
of $(Y, q)$ and vertices correspond to $p_0$. 
The following remark is important.

\rem{\label{simple_cylinder_principle}
In the construction of  triangles,  
as we do in Section \ref{proof_of_main_thm}, 
the last triangle has two sides 
corresponding to returning geodesics of $(Y, q)$.
 %Otherwise, we can proceed the construction.
 The preimage of the triangle via $\varphi$ 
 is a simple cylinder that is invariant under $\tau$.
 }

Moreover, we use the following lemma.

\lem{\label{end_lem}
Let $(Y_0, q_0)$ be a flat surface of genus $0$ 
such that $q_0$ has no zero and $4$ simple poles $p_0, p_1, p_2, p_3$.
Let  $s$ be a saddle connection on $(Y, q)$ connecting $p_0$ and $p_1$.
Let $\vv$ be a unit vector that is not parallel to $s$. 
Then there exists  
a (closed) triangle $\Delta$ in $\mathbb{C}$ 
with sides $s_0, s_1, s_2$ in counterclockwise order 
and an  orientation preserving  immersion 
$\rho : \Delta \to (Y, q)$ 
 satisfying the following: 
\begin{enumerate}
\item $s_0$ is parallel to $s$ and $\Delta$ is a left (resp. right) strongly $(s_0, \vv)$-restricted triangle, 
\item $\rho|_{{\rm Int}(\Delta)}$ is a local isometric embedding,
\item $\rho({\rm Int}(\Delta))$ contains no singular points, 
\item every vertex of $\Delta$ is mapped to $p_0$,
\item $s=\rho \circ s_0$, and 
\item the preimage of $\left\{ p_2, p_3 \right\}$ 
coincides with the set of midpoints of $s_1$ and $s_2$.
%\item If $s_i$ is parallel to $\vv$ for some $i=1,2$, then $\rho \circ s_i$ is a returning geodesic.
\end{enumerate}
}

\pf{
We prove in the case where $\Delta$ is  
a left strongly $(s_0, \vv)$-restricted triangle.
%To make the argument easy, we assume that $s$ is horizontal.
Cut $(Y_0, q_0)$ along $s$.
Then $(Y_0,q_0)$ 
is a flat surface with a geodesic boundary $s$.
There exists an Euclidean cylinder $C$ in $(Y_0, q_0)$ 
one of whose boundary components is this geodesic boundary $s$.
We can extend $C$ so that the other boundary component 
%of $C$ than $s$ 
contains $p_2$ or $p_3$. 
Assume that this boundary component $C$ contains only $p_2$. 
Then $Y_0 - C$ is not empty. 
However, it contradicts the fact that the angle around $p_2$ is $\pi$. 
Therefore, this boundary component contains both $p_2$ and $p_3$. 
Moreover, $Y_0$ is obtained from $C$ by sewing each boundary component.
We cut $C$ along the segment from $p_0$ that is parallel to $\vv$.
Then the resulting surface $P$ is regarded as a parallelogram in $\mathbb{C}$ 
with sides that are parallel to $s_0$ 
and there exists a natural map $\rho_0 : P \to (Y_0, q)$. 
Let $s$ be a side of $P$ that is mapped  to $s_0$ via $\rho_0$. 
Let $a$ be the side of $P$ that is opposite to $s$. 
We label the other side of $P$ with $b$ and $c$ 
so that the counterclockwise  order of sides is $s, c, a, b$.
Let $P_0$ be the vertex of $P$ 
that is the intersection of $s$ and $b$ 
and
$P_1$ the vertex of $P$ 
that is the intersection of $s$ and $c$.  
Then $\rho_0 \circ a$ is a saddle connection 
connecting $p_2$ and $p_3$. 
The sets $\rho_0^{-1}(p_2)-b$ 
and $\rho_0^{-1}(p_3)-b$
contain only one point. 
We set $\rho_0^{-1}(p_i)-b=\left\{P_i \right\}$ 
for $i=2, 3$.
We may assume without loss of generality that 
$P_0P_2 <P_0P_3$.
Let $P_4$ be the intersection of the lines $P_0P_2$ 
and $P_1P_3$. 
Then $\Delta=\triangle P_0P_1P_4$ is a
 left strongly $(s_0, \vv)$-restricted triangle.
 We set $s_1=P_1P_4$ and $s_2=P_0P_4$.  
Let $l$ be the line 
that  passes through $P_4$ 
and is parallel to $\vv$ 
(see Figure \ref{end_lem_pic}).
If $P_3 \not \in c$, then 
the line $l$ divides $\Delta-\textrm{Int}(P)$ into two 
components $\Delta_1$ and $\Delta_2$. 
We assume that $\Delta_1$ contains $P_2$.
Then the map $\rho : \Delta \to (Y_0, q_0)$ defined by 
 $$\rho(z)=
\begin{cases}
\rho_0 (z)   & \left(z \in P \right) \\
\rho_0 \circ r_{P_2}(z) & \left(z\in \Delta_1 \right)\\
\rho_0 \circ r_{P_3} (z) & \left(z\in \Delta_2 \right).
\end{cases}
$$
satisfies all conditions as in the claim.
Here, $r_{P_2}$ and $r_{P_3}$ 
are point reflections in $P_2$ and $P_3$, respectively.
If $P_3 \in c$, then $P_2$ is the midpoint of $a$ 
and $l$ coincides with the line $P_1P_3$. 
Then the map $\rho : \Delta \to (Y_0, q_0)$ defined by 
 $$\rho(z)=
\begin{cases}
\rho_0 (z)   & \left(z \in P \right) \\
\rho_0 \circ r_{P_2}(z) & \left(z\in \Delta-P \right)
\end{cases}
$$
satisfies all conditions as in the claim.
%Let $s^\prime $ be a segment in $C$ 
%connecting $p_0$ to $p_2$ 
%There exists a segment $s^{\prime\prime}$ in $C$ 
%that 
%connects$p_0$ to $p_3$ 
%and does not intersects with $s^\prime$ other than $p_0$.
%Cut $C$ along $s^\prime$ and $s^{\prime\prime}$, 
%$C$ is divided into two components.
%Glue them along the sides containing $p_2$ and $p_3$, 
%we obtain a triangle $\Delta^\prime$ that can immersed to $(Y_0, q_0)$.
%Let $s_0$ be the side of $\Delta^\prime$ corresponding to $s$ and 
%set $A=
%\left[
%\begin{array}{cc} 
%1& |s| \\ 
%0&1 
%\end{array}
%\right]
%$. 
%Then, there exists $n\in \mathbb{Z}$ such that 
%$\Delta=A^n\left (\Delta^\prime \right)$
%is a left strongly $(s_0, \vv)$-restricted triangle. 
%Let denote $s_1$ and $s_2$ the sides of $\Delta$ edaLmargin=0 $s_0$ 
%so that $s_0, s_1$ and $s_2$ are in counterclockwise order.
%By construction, we can obtain  an  orientation preserving  immersion 
%$\rho : \Delta \to (Y, q)$ 
%satisfying all conditions we want.
\begin{figure}[ht!]
\labellist
%\small
\hair 0pt
%\pinlabel $O$  at 210 -15
\pinlabel $s$  at 250 0
\pinlabel $a$  at 285 170
\pinlabel $b$  at 40 115
\pinlabel $c$  at 500 115
\pinlabel $\vv$  at 10 60
\pinlabel $P_0$  at -20 15
\pinlabel $P_1$  at 435 15
\pinlabel $P_2$  at 145 175
\pinlabel $P_3$  at 380 175
\pinlabel $P_4$  at 300  300
\pinlabel $l$  at 340 315
 \endlabellist
\centering
 \includegraphics[scale=0.4]{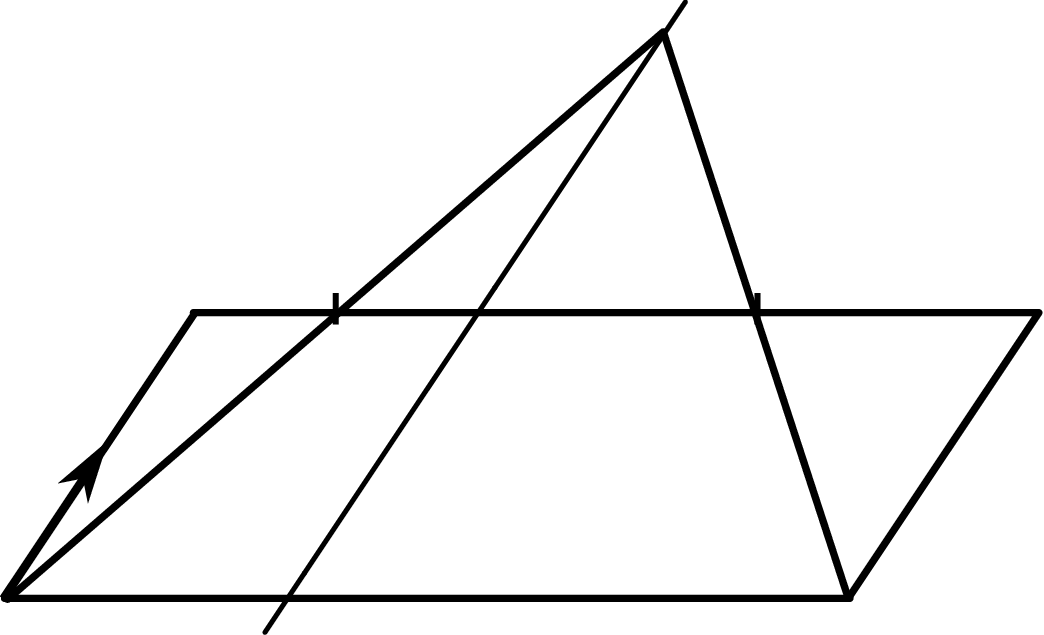} 
\caption{The sides $s, a, b, c $, the points $P_0, P_1, P_2, P_3, P_4$,  and the line $l$.}
\label{end_lem_pic}
\end{figure}
}

\begin{figure}[h]
\labellist
\hair 0pt
\pinlabel $s_0$  at   0 40
\pinlabel $s_1$  at   110 40
\pinlabel $s_2$  at   80 105
\pinlabel $s_3$  at   170 180
\pinlabel $s_4$  at   60 180
\pinlabel $\Delta_0$  at   50 60 
\pinlabel $\Delta_1$  at 115 155
\endlabellist
\centering
\includegraphics[scale=0.38]{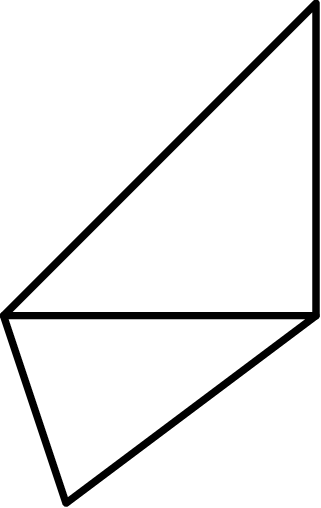} 
\caption{The triangles $\Delta_0$ and $\Delta_1$.}
\label{classification01}
\end{figure}

\begin{figure}[h]
\labellist
\hair 0pt
\pinlabel $s_0$  at   0 40
\pinlabel $s_1$  at   110 40
\pinlabel $s_2$  at   80 105
\pinlabel $s_3$  at   170 180
\pinlabel $s_4$  at   60 180
\pinlabel $\Delta_0$  at   50 60 
\pinlabel $\Delta_1$  at 115 155
\pinlabel $\Delta_2$  at  200 210
\pinlabel $s_5$  at   240 180
\pinlabel $s_6$  at   205 255
\endlabellist
\centering
\includegraphics[scale=0.38]{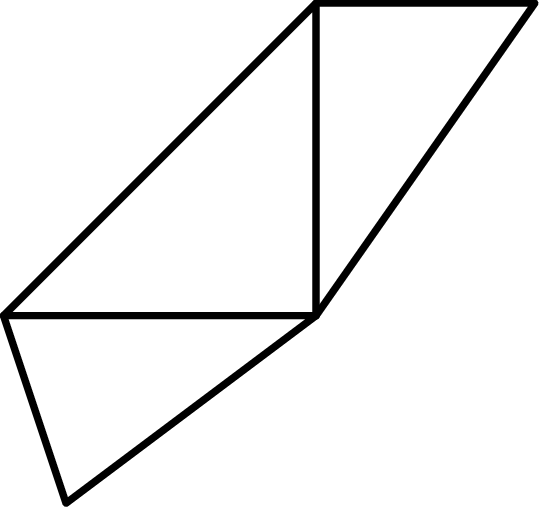} 
\caption{The triangles $\Delta_0$,  $\Delta_1$, and $\Delta_2$.}
\label{classification012}
\end{figure}

\begin{figure}[h]
\labellist
\hair 0pt
\pinlabel $s_0$  at   0 40
\pinlabel $s_1$  at   110 40
\pinlabel $s_2$  at   80 105
\pinlabel $s_3$  at   170 180
\pinlabel $s_4$  at   60 180
\pinlabel $\Delta_0$  at   50 60 
\pinlabel $\Delta_1$  at 115 155
\pinlabel $\Delta_2$  at  200 210
\pinlabel $s_5$  at   240 180
\pinlabel $s_6$  at   205 255
\pinlabel $s_7$  at   280 325
\pinlabel $s_8$  at   175 325
\pinlabel $\Delta_3$  at  230 315
\endlabellist
\centering
\includegraphics[scale=0.38]{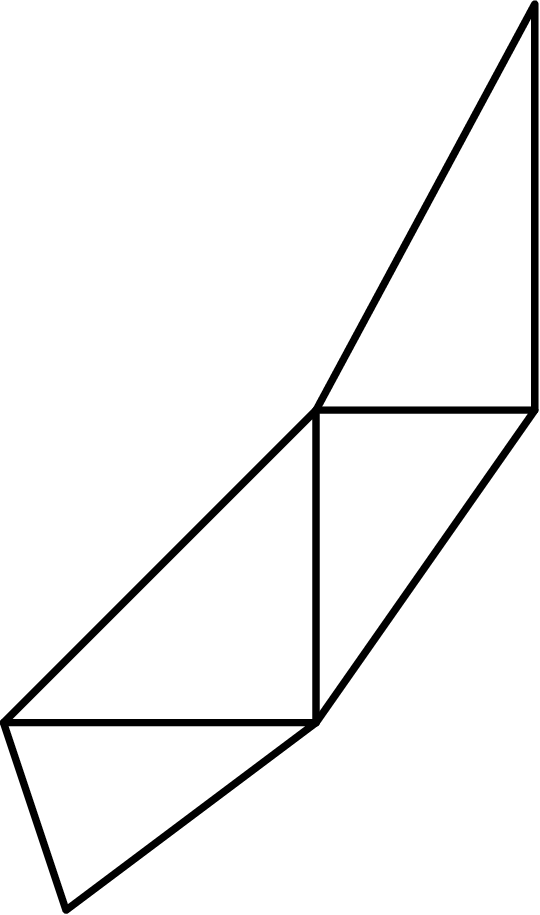} 
\caption{The triangles $\Delta_0$,  $\Delta_1$, $\Delta_2$, and $\Delta_3$.}
\label{classification0123}
\end{figure}

\begin{figure}[h]
\labellist
\hair 0pt
\pinlabel $s_0$  at   -20 90
\pinlabel $s_1$  at   120 90
\pinlabel $s_2$  at   80 165
\pinlabel $s_3$  at   170 240
\pinlabel $s_4$  at   60 240
\pinlabel $\Delta_0$  at   50 110 
\pinlabel $\Delta_1$  at 115 215 %115 155
\pinlabel $\Delta_2$  at  200 270
\pinlabel $s_5$  at   270 240
\pinlabel $s_6$  at   225 315
\pinlabel $s_7$  at   325 385
\pinlabel $s_8$  at   205 385
\pinlabel $\Delta_3$  at  270 375
\pinlabel $s_9$  at   410 385
\pinlabel $s_{10}$  at  390 470
\pinlabel $\Delta_4$  at 370 420
\endlabellist
\centering
\includegraphics[scale=0.38]{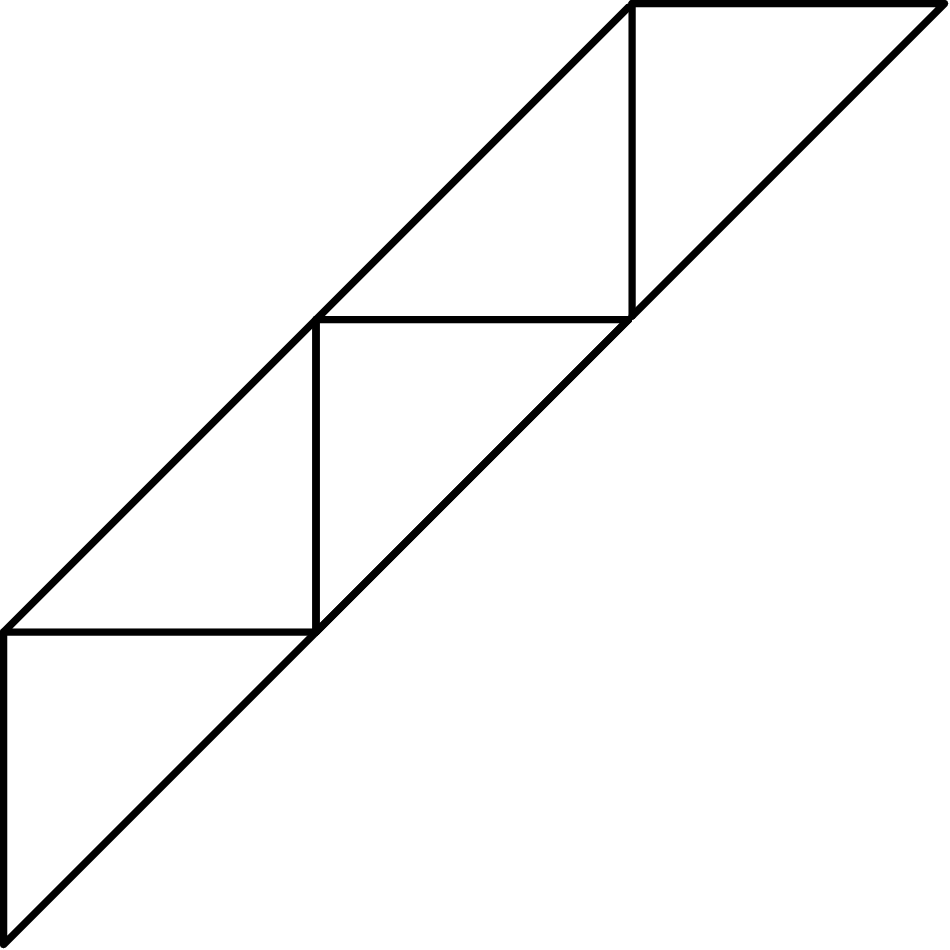} 
\caption{The triangles $\Delta_0$,  $\Delta_1$, $\Delta_2$, $\Delta_3$, and $\Delta_4$ that construct $(Y, q)$.}
\label{classification01243}
\end{figure}

For the proof of Theorem \ref{classification}, we will construct triangles 
$\Delta_0, \Delta_1, \Delta_2, \dots, $ 
and their sides 
$s_0, s_1, s_2, \dots$.
Hereafter, 
$Q_i$ denotes the midpoint of the side $s_i$ 
for $i=1,2, \dots$. 
Given a returning geodesic $s$ on $(Y, q)$. 
By Theorem \ref{good_triangle}, 
there exists a triangle $\Delta_0$ in $\mathbb{C}$ 
with sides $s_0, s_1, s_2$ in counterclockwise order 
and an orientation preserving  immersion 
$\rho_0 : \Delta_0 \to (Y, q)$ 
 satisfying all conditions 
 as in the claim of Theorem \ref{good_triangle}.
 The preimage $C_1=\varphi^{-1}(\rho_0(\Delta_0))$ 
 is a simple cylinder 
 that is invariant under $\tau$.
The surface $Y_1=Y-\rho_0(\Delta_0)$ 
is of genus $0$ and has a boundary $\rho_0\circ s_2$. 
%Let $q_1$ be the restriction of $q$ onto $Y_1$. 
Sewing the boundary, 
$Y_1$ is a closed surface of genus $0$  
and $q$ induces  a meromorphic quadratic differential 
$q_1$ on $Y_1$ such that 
 $\rho\circ s_2$ is a returning geodesic 
on the flat surface $(Y_1, q_1)$. 
Applying Theorem \ref{good_triangle}, 
there exists a triangle $\Delta_1$ in $\mathbb{C}$ 
with sides $s_2, s_3, s_4$ in counterclockwise order 
and an orientation preserving  immersion 
$\rho_1 : \Delta_1 \to (Y_1, q_1)$ 
 satisfying the following;
\begin{enumerate}[(\textrm{1-}1)]
\item $\rho_1|_{\interior{\Delta_1}}$ is a local isometric embedding,
\item $\Delta_1 \cap \Delta_0 = s_2$,
\item $\rho_1|_{\interior{\Delta_1}}$ contains no singular points, 
\item every vertex of $\Delta_1$ is mapped to $p_0$,
\item $\rho_1 \circ s_2 =\rho_0 \circ s_2$, and 
\item\label{can_replace} $\rho_1 \circ s_4$ is a returning geodesic.
\end{enumerate}
Rename $A\cdot (Y, q)$ 
to $(Y, q)$ 
for some $A\in \gltr$ if necessary, 
we may assume that 
$s_2$ is horizontal, 
$s_3$ is vertical,  
$|s_2|=|s_3|$, 
and $\Delta_1$ is above $\Delta_0$ (see Figure \ref{classification01}).
The surface $Y_2=Y_1-\rho_1(\Delta_1)$ 
is of genus $0$ and has a boundary $\rho_1 \circ s_3$. 
Sewing the boundary, 
$Y_2$ is a closed surface of genus $0$  
and $q_1$ induces  a meromorphic quadratic differential 
$q_2$ on $Y_2$ such that 
 $\rho_1 \circ s_3$ is a returning geodesic 
on the flat surface $(Y_2, q_2)$. 
Applying Theorem \ref{triangle} 
with the vector $\vv=(1, 0)$,  
there exists  
a  triangle $\Delta_2 $ in $\mathbb{C}$ 
with sides $s_3, s_5, s_6$ in counterclockwise order 
and an  orientation preserving  immersion 
$\rho_2 : \Delta_2 \to (Y_2, q_2)$ 
 satisfying the following: 
\begin{enumerate}[(\textrm{2-}1)]
\item 
%$s_$ is parallel to $s$ and 
$\Delta_2$ is a right strongly $(s_3, \vv)$-restricted triangle, 
\item $\rho_2|_{{\rm Int}(\Delta_2)}$ is a local isometric embedding,
\item $\Delta_2 \cap \Delta_1 = s_3$,
\item $\rho_2({\rm Int}(\Delta_2))$ contains no singular points, 
\item every vertex of $\Delta_2$ is mapped to $p_0$,
\item $\rho_2 \circ s_3=\rho_1 \circ s_3$, and 
\item $\rho_2 \circ s_i$ is a saddle connection or a returning geodesic for each $i=5, 6$.
%\item If $s_i$ is parallel to $\vv$ for some $i=1,2$, then $\rho \circ s_i$ is a returning geodesic.
\end{enumerate}
If both $\rho_2 \circ s_5$ and $\rho_2 \circ s_6$ 
are saddle connections, 
$Y_2-\rho_2 (\Delta_2)$ 
has two connected components. 
By Remark \ref{simple_cylinder_principle}, 
we can find 
triangles whose preimages via $\varphi$ 
are simple cylinders 
that are invariant under $\tau$
in the connected components.
This contradicts the assumption that $(Y, q)$ has 
 no three
disjoint simple cylinders each of which  is 
invariant under $\tau $. 
Since $q$ has $7$ simple poles, 
only one of 
$\rho_2 \circ s_5$ and  $\rho_2 \circ s_6$ 
is a saddle connection. 
Moreover, by Remark \ref{simple_cylinder_principle}, 
we can find a triangle  in $Y_2-\rho_2(\Delta_2)$
whose preimage via $\varphi$ 
is a simple  cylinder 
that is invariant under $\tau$.
If $\rho_2 \circ s_6$ is a returning geodesic, 
the segment $Q_4Q_6$  
is mapped via $\rho_1$ and $\rho_2$ 
to a saddle connection  in $(Y, q)$ 
whose preimage via $\varphi$ is a core curve of a simple cylinder that is invariant under $\tau$.
This contradicts the assumption that $(Y, q)$ has 
 no three
disjoint simple cylinders each of which  is 
invariant under $\tau $. 
Thus,  $\rho_2 \circ s_5$ is a returning geodesic.
The segment $Q_4Q_5$ cannot correspond 
to a saddle connection in $(Y, q)$ 
whose preimage via $\varphi$ is a core curve of a simple cylinder. 
This implies that $Q_4Q_5$ and $s_6$ are horizontal 
(see Figure \ref{classification012}).
The surface $Y_3=Y_2-\rho_2(\Delta_2)$ 
is of genus $0$ and has a boundary $\rho_2 \circ s_6$. 
Sewing the boundary, 
$Y_3$ is a closed surface of genus $0$  
and $q_2$ induces  a meromorphic quadratic differential 
$q_3$ on $Y_3$ such that 
 $\rho_2 \circ s_6$ is a returning geodesic 
on the flat surface $(Y_3, q_3)$. 
Applying Theorem \ref{triangle} 
with the vector $\vv_1=(0, 1)$,  
there exists  
a  triangle $\Delta_3 $ in $\mathbb{C}$ 
with sides $s_6, s_7, s_8$ in counterclockwise order 
and an  orientation preserving  immersion 
$\rho_3 : \Delta_3 \to (Y_3, q_3)$ 
 satisfying the following: 
\begin{enumerate}[(\textrm{3-}1)]
\item 
%$s_$ is parallel to $s$ and 
$\Delta_3$ is a left strongly $(s_6, \vv_1)$-restricted triangle, 
\item $\rho_3|_{{\rm Int}(\Delta_3)}$ is a local isometric embedding,
\item $\Delta_3 \cap \Delta_2 = s_6$,
\item $\rho_3({\rm Int}(\Delta_3))$ contains no singular points, 
\item every vertex of $\Delta_3$ is mapped to $p_0$,
\item $\rho_3 \circ s_6=\rho_2 \circ s_6$, and 
\item $\rho_3 \circ s_i$ is a saddle connection or a returning geodesic for each $i=7, 8$.
%\item If $s_i$ is parallel to $\vv$ for some $i=1,2$, then $\rho \circ s_i$ is a returning geodesic.
\end{enumerate}
By the same argument as above, 
we see that 
$\rho_3 \circ s_8$ is a returning geodesic, 
$\rho_3 \circ s_7$ is a vertical saddle connection, and
$s_7$ is vertical (see Figure \ref{classification0123}). 
Assume that $|s_6|<|s_2|$. 
Let $\vv_2$ 
be  the unit vector 
that is parallel to $s_4$ 
from top to bottom. 
Then we reconstruct $\Delta_0$ 
so that $\Delta_0$ is a 
right strongly $(s_2, \vv_2)$-restricted triangle.
Then $Q_0 Q_4$ and $Q_1 Q_5$ 
are mapped via $\rho_0$, $\rho_1$, and $\rho_2$ 
to saddle connections  in $(Y, q)$ 
whose preimages via $\varphi$ are core curves of  
simple cylinders that are invariant under $\tau$.
Moreover, by Remark \ref{simple_cylinder_principle}, 
we can find a triangle  in $Y_2-\rho_2(\Delta_2)$
whose preimage via $\varphi$ 
is a simple cylinder 
that is invariant under $\tau$.
This contradicts the assumption that $(Y, q)$ has 
 no three
disjoint simple cylinders each of which  is 
invariant under $\tau $. 
Thus, we have $|s_6| \geq |s_2|$.
Next, if  $|s_7|<|s_3|$, 
then the segment $Q_4Q_8$ 
is mapped via $\rho_1$, $\rho_2$, and $\rho_3$ 
to a saddle connection  in $(Y, q)$ 
whose preimage via $\varphi$ is a core curve of a simple cylinder that is invariant under $\tau$.
Moreover, by Remark \ref{simple_cylinder_principle}, 
we can find a triangle  in $Y_3-\rho_3(\Delta_3)$
whose preimage via $\varphi$ 
is a simple cylinder 
that is invariant under $\tau$.
This contradicts the assumption that $(Y, q)$ has 
 no three
disjoint simple cylinders each of which  is 
invariant under $\tau $.  
Thus, we have $|s_7| \geq |s_3|$.
The surface $Y_4=Y_3-\rho_3(\Delta_3)$ 
is of genus $0$ and has a boundary $\rho_3 \circ s_7$. 
Sewing the boundary of $Y_4$, 
$Y_4$ is a closed surface of genus $0$  
and $q_3$ induces  a meromorphic quadratic differential 
$q_4$ on $Y_4$ such that 
 $\rho_3 \circ s_7$ is a returning geodesic 
on the flat surface $(Y_4, q_4)$. 
Let $\vv_3$ be the unit vector 
that is parallel to $s_5$ 
from bottom to top. 
Applying Lemma \ref{end_lem} 
with the vector $\vv_3$,  
there exists  
a  triangle $\Delta_4^\prime $ in $\mathbb{C}$ 
with sides $s_7, s_9^\prime, s_{10}^\prime$ in counterclockwise order 
and an  orientation preserving  immersion 
$\rho_4^\prime : \Delta_4^\prime \to (Y_4, q_4)$ 
 satisfying the following: 
\begin{enumerate}[(\textrm{4-}1)$^\prime$]
\item 
%$s_$ is parallel to $s$ and 
$\Delta_4^\prime$ is a right strongly $(s_7, \vv_3)$-restricted triangle, 
\item $\rho_4^\prime|_{{\rm Int}(\Delta_4^\prime )}$ is a local isometric embedding,
\item $\Delta_4^\prime \cap \Delta_3 = s_7$,
\item $\rho_4^\prime ({\rm Int}(\Delta_4^\prime))$ contains no singular points, 
\item every vertex of $\Delta_4^\prime $ is mapped to $p_0$,
\item $\rho_4^\prime \circ s_7=\rho_3 \circ s_7$, and 
\item $\rho_4^\prime \circ s_i^\prime$ is a returning geodesic for each $i=9, 10$.
%\item If $s_i$ is parallel to $\vv$ for some $i=1,2$, then $\rho \circ s_i$ is a returning geodesic.
\end{enumerate}
Let $Q_{i}^\prime$ be the midpoint if $s_i^\prime$ for $i=9, 10$. 
Then $Q_5 Q_9^\prime$ 
is mapped via $\rho_1$, $\rho_2$, $\rho_3$, and $\rho_4^\prime$  
to a saddle connection  in $(Y, q)$ 
whose preimage via $\varphi$ is a core curve of a simple cylinder that is invariant under $\tau$.
By  the assumption that 
$(Y, q)$ has 
 no three
disjoint simple cylinders each of which  is 
invariant under $\tau $,  
$s_8$ must be parallel to $s_{10}^\prime$.  
If not, 
$Q_8 Q_{10}^\prime$ 
is mapped via 
%$\rho_1$, $\rho_2$, 
$\rho_3$ and $\rho_4^\prime$  
to a saddle connection  in $(Y, q)$ 
whose preimage via $\varphi$ is a core curve of a simple cylinder that is invariant under $\tau$. 
Therefore, we have  $|s_7| = |s_3|$. 
Moreover, if $|s_6|>|s_2|$, 
then the segment $Q_4Q_8$ 
is mapped via $\rho_1$, $\rho_2$, and $\rho_3$ 
to a saddle connection  in $(Y, q)$ 
whose preimage via $\varphi$ is a core curve of a simple cylinder that is invariant under $\tau$.
This contradicts the assumption that $(Y, q)$ has 
 no three
disjoint simple cylinders each of which  is 
invariant under $\tau $.  
Thus, we have $|s_6|=|s_2|$.
Again, applying Lemma \ref{end_lem} 
with the vector $\vv_3$,  
there exists  
a  triangle $\Delta_4 $ in $\mathbb{C}$ 
with sides $s_7, s_9, s_{10}$ in counterclockwise order 
and an  orientation preserving  immersion 
$\rho_4 : \Delta_4 \to (Y_4, q_4)$ 
 satisfying the following: 
\begin{enumerate}[(\textrm{4-}1)]
\item 
%$s_$ is parallel to $s$ and 
$\Delta_4$ is a left strongly $(s_7, \vv_3)$-restricted triangle, 
\item $\rho_4|_{{\rm Int}(\Delta_4)}$ is a local isometric embedding,
\item $\Delta_4 \cap \Delta_3 = s_7$,
\item $\rho_4({\rm Int}(\Delta_4))$ contains no singular points, 
\item every vertex of $\Delta_4$ is mapped to $p_0$,
\item $\rho_4 \circ s_7=\rho_3 \circ s_7$, and 
\item $\rho_4 \circ s_i$ is a returning geodesic for each $i=9, 10$.
%\item If $s_i$ is parallel to $\vv$ for some $i=1,2$, then $\rho \circ s_i$ is a returning geodesic.
\end{enumerate}
Then $Q_8Q_{10}$ 
is mapped via 
%$\rho_1$, $\rho_2$, 
$\rho_3$ and $\rho_4$ 
to a saddle connection  in $(Y, q)$ 
whose preimage via $\varphi$ is a core curve of a simple cylinder that is invariant under $\tau$.
Thus, the side $s_9$ must be parallel to $s_5$. 
Moreover, the segment $Q_4 Q_9$ 
must not correspond to a core curve of a simple cylinder that is invariant under $\tau$.
This implies that $s_{10}$ is horizontal.
Now, $\Delta_1$, $\Delta_2$, $\Delta_3$, and $\Delta_4$ 
are all equilateral right triangles. 
By the same argument as above, 
we can reconstruct $\Delta_0$ 
to be an equilateral right triangles equilateral right triangles 
with a vertical side $s_0$. 
Therefore, $(Y, q)$ is the flat surface 
constructed from $5$  right triangles equilateral right triangles 
as in Figure \ref{classification01243}  
so that 
the midpoints of the sides $s_0$, $s_1$, $s_4$, $s_5$, $s_8$, $s_9$, and $s_{10}$ 
correspond to simple poles of $(Y, q)$.  
This implies that 
$(X, \omega)$ coincides with $\str$.

%\input{Appendix}

%%%%%%%%%%%%%%%%%%%%   End of main body of article
%
%                             References
%
%   BiBTeX users uncomment the following line:
%
\bibliographystyle{alpha}
%

%

%\begin{thebibliography}
\bibliography{ref}
%\end{thebibliography}

\end{document}